\documentclass[12pt]{article}

\usepackage{amsfonts}
\usepackage{graphics}
\usepackage{amsmath}
\usepackage{amscd}
\usepackage{rlepsf}
\usepackage{color}

\usepackage{euscript}

\newtheorem{Theorem}{Theorem}

\newtheorem{Corollary}[Theorem]{Corollary}
\newtheorem{Example}[Theorem]{Example}
\newtheorem{Definition}[Theorem]{Definition}

\newtheorem{Remark}[Theorem]{Remark}

\newtheorem{Lemma}[Theorem]{Lemma}
\newtheorem{Proposition}[Theorem]{Proposition}

\newtheorem{Fundamental Theorem}{Fundamental Theorem}

\newenvironment{Proof}[1][Proof]{\textbf{#1.} }{\ \rule{0.5em}{0.5em}}

\def \A {{\cal A}}
\def \e {\epsilon}

\def \i {\iota}

\def \l {{\leq t}}
\def \n {\nu}

\def \Hom {\mathrm{Hom}}
\def \Z {\mathbb{Z}}
\def \id {\mathrm{id}}

\def \f {\phi}
\def \p {\psi}

\def \ra {\xrightarrow}

\def \g {\gamma}
\def \G {\mathcal{G}}
\def \d {\partial}
\def \t {\triangleright}

\def \S {\Sigma}
\def \R {\mathbb{R}}
\def \s {\scriptstyle}

\def \N {\mathbb{N}_0}
\begin{document}

\title{{The Fundamental Crossed Module of the Complement of a Knotted Surface}}

\author{Jo\~{a}o  Faria Martins \footnote{Also at: Universidade Lus\'{o}fona de Humanidades e Tecnologia, Av do Campo Grande, 376, 1749-024, Lisboa, Portugal. }\\ \footnotesize\it  {Departamento de Matem\'{a}tica, Instituto Superior T\'{e}cnico,}\\ {\footnotesize\it Av. Rovisco Pais, 1049-001 Lisboa, Portugal}\\ {\footnotesize\it jmartins@math.ist.utl.pt}}

\date{\today}

\maketitle
\begin{abstract}
We prove that if $M$ is a CW-complex and  $M^1$ is its 1-skeleton then the crossed module $\Pi_2(M,M^1)$  depends  only on the homotopy type of $M$ as a space,  up to free products{, in the category of crossed modules,} with $\Pi_2(D^2,S^1)$. From this it follows that,  if $\G$ is a finite crossed module and  $M$ is finite, then the number of crossed module morphisms  $\Pi_2(M,M^1) \to \G$ can be re-scaled to a homotopy invariant $I_\G(M)$, depending only on the homotopy  2-type of $M$.  We describe an  algorithm {for calculating} $\pi_2(M,M^{(1)})$ as a crossed module over
  $\pi_1(M^{(1)})$, in the case when $M$ is the complement of a knotted surface $\S$
  in $S^4$ and $M^{(1)}$ is the handlebody made from {the} $0$- and $1$-handles of a handle decomposition of $M$. {Here $\S$ is presented by a knot with bands}. This  in particular {gives us  a  geometric method {for calculating} the algebraic 2-type of the complement of a knotted surface from a hyperbolic splitting of it.} We prove in addition  that  the invariant $I_\G$ yields a non-trivial invariant of knotted surfaces in $S^4$ with good properties with {regards} to explicit calculations. \end{abstract}
\tableofcontents
\section*{Introduction}
Let $(M,N)$ be a pair of based path-connected spaces.
The concept of a crossed module arises from  a universal description of the properties of the boundary map  $\d\colon \pi_2(M,N) \to \pi_1(N)$, together   with the natural action of $\pi_1(N)$ on $\pi_2(M,N)$. These  data define the crossed module $\Pi_2(M,N)$, called {the} ``Fundamental Crossed Module of $(M,N)$''.

 Due to some strong theorems by  J.H.C. Whitehead, it is possible, in principle, to calculate $\Pi_2(M,M^1)$, when  $M$ is a connected CW-complex and $M^1$ is its 1-skeleton.  The calculability of fundamental crossed modules is, in addition, enhanced by a 2-dimensional van Kampen theorem due to R. Brown and P.J. Higgins, stating that, under mild conditions,   the fundamental crossed module functor  {from the category of based pairs of path connected spaces to the category of crossed modules} preserves colimits; see \cite{BH,BH3,B2}.

{The crossed module} $\Pi_2(M,M^1)$ determines not only  $\pi_1(M)$ and $\pi_2(M)$ as a module over $\pi_1(M)$,  but also it determines the $k$-invariant $k(M)\in H^3(\pi_1(M),\pi_2(M))$; in other words all the algebraic 2-type of $M$, thus {it is} strictly stronger than $\pi_1(M)$ and $\pi_2(M)$ alone.

Crossed modules admit an obvious notion of homotopy; see {\cite{W3,BH2,FMP}.} From the point of view of the homotopy type of  a CW-complex $M$, as a space, considering {its fundamental crossed module $\Pi_2(M,M^1)$}  introduces an ambiguity {due to} the choice of a cell decomposition. However, despite this asymmetry, the homotopy type of the fundamental crossed module $\Pi_2(M,M^1)$ depends only on the homotopy type of $M$ as a space; see \cite{W3,BH2}. In fact this result can be improved. {It is proved in this article that  the crossed module $\Pi_2(M,M^1)$} does not  depend on the CW-decomposition of $M$, up to free products{, in the category of crossed modules,} with  $\Pi_2(D^2,S^1)$. {A similar result can be obtained for the fundamental crossed complex of a CW-complex; see \cite{FM2}.} {Consequently,} the fundamental crossed module $\Pi_2(M,M^1)$, up to free {products} with $\Pi_2(D^2,S^1)$, is a genuine homotopy invariant of CW-complexes. This invariant determines the topological  2-type of $M$; see \cite{Bau,Lo}. 

{In the light of the previous discussion, it is} natural to consider crossed modules in Geometric Topology.  In this article we present an algorithm {for calculating} the crossed module  $\Pi_2(M,M^{(1)})$ in the case when $M$ is the complement of a knotted surface $\S$ in $S^4$ and $M^{(1)}$ is the 1-handlebody {(made out of the $0$- and $1$-handles)} of a handle decomposition of {$M$}. Both the handle decomposition of $M$ (following \cite{G,GS,CKS})  and the method  {for  determining} $\Pi_2(M,M^{(1)})$ are defined  from a hyperbolic splitting  of $\S$.  We thus have a completely geometric {method for calculating} the algebraic 2-type of the complement of a knotted surface from a hyperbolic splitting of it.  {This will be the main result of this article.}  {Our method to determine the algebraic 2-type of the complement of a knotted surface {in $S^4$} should be compared with Lomonaco's method in \cite{L}. See also \cite{M}.}

{We illustrate our {algorithm}  by {calculating  $\Pi_2(M,M^{(1)})$}, where $M$ is the complement of a knotted surface, in some non-trivial cases. We use it {for determining} {the second homotopy group of  the Spun Trefoil complement}, and verify {that the result} agrees with previous calculations.}

Given that the fundamental crossed module $\Pi_2(M,M^1)$ {depends only} on the homotopy type of $M$, as a space, up to free products with  $\Pi_2(D^2,S^1)$,  it follows that, if $\G$ is a finite crossed module and $M$ is a finite CW-complex, then the number of crossed module morphisms {$\Pi_2(M,M^1) \to \G$} (which is finite) can be re-scaled to a homotopy invariant $I_\G(M)$ (the ``Crossed Module Invariant''), thus solving again the problem that we introduced an ambiguity {due to} the choice of {a CW-decomposition of $M$.} {This {Crossed Module Invariant} depends only on the topological 2-type of a space.} {This gives an extension of Yetter's invariant of 3-manifolds \cite{Y} to general CW-complexes. This invariant can also be extended to depend on a crossed complex (more general than a crossed module); see \cite{FM2,FMP}.}

From the algorithm {for calculating} $\Pi_2(M,M^{(1)})$ when $M$ is the complement of an embedded surface in $S^4${, a method follows for calculating} $I_\G(M)$, where $\G$ is a finite crossed module.  We  prove that the Crossed Module Invariant  $I_\G$ defines a non trivial invariant of knotted surfaces, with good calculability properties, especially in the case of abelian crossed modules. {We will give an example (appearing also in \cite{FMK}) of a pair of knotted surfaces, each diffeomorphic to the disjoint union of two {tori} $S^1 \times S^1$, with the same fundamental group of the complement, but distinguished by their {crossed module invariants}.  An open problem is whether the {Crossed Module Invariant} is strong enough to distinguished between knotted surfaces with the same fundamental {and} second homotopy {groups of the complement (seen as  $\pi_1$-modules)}, but with distinct $k$-invariants. See also \cite{PS}.}

In \cite{FM1} we defined an invariant of knotted surfaces from any  finite crossed module. The
construction was inspired by previous work of D. Yetter and T. Porter on
manifold invariants defined from models of homotopy $2$-types (crossed modules
of groups), see \cite{Y} and \cite{P,P2}. This
article should give, in particular, an interpretation of our previous construction.

 \section{Preliminaries and General Results}

\subsection{Crossed Modules}\label{general}
Let $G$ and $E$ be groups. {A  crossed module $\G=(G,E,\d,\t)$} is given by a group morphism $\d \colon  E \to G$ and an action $\t$ of $G$ on $E$ on the left by automorphisms. The conditions on $\t$ and $\d$ are:
\begin{enumerate}
\item $\d(X \t e)=X\d(e)X^{-1};\forall X \in G, \forall e \in E$,
\item $\d(e) \t f=e f e^{-1}; \forall e, f \in E$;
\end{enumerate}
see for example  \cite{B2,BV}.
Notice that the second condition implies that ${\rm ker} (\d)$ commutes with all $E$. Therefore, the action of $G$ on $E$ induces an action of ${\rm coker}(\d)$ on ${\rm ker} (\d)$, by automorphisms.
We call $G$ the base group and $E$ the principal group. A crossed module is
called finite if both $G$ and $E$ are finite. 

The significance of the definition of crossed modules for Geometric Topology stems  from: 
\begin{Example}\label{crs}
  Let $(M,N)$ be a pair of based {path connected spaces}. Then  $\left (\pi_1(N),\pi_2(M,N), \d,\t \right
 )\doteq\Pi_2(M,N)$, where the boundary map {$\d \colon \pi_2(M,N)\to \pi_1(N)$ and the action of $\pi_1(N)$ on $\pi_2(M,N)$ are the natural} ones,  is a crossed module, called the ``Fundamental Crossed Module'' of $(M,N)$.  This is a result of J.H.C. Whitehead, see \cite{W1}.  
\end{Example}
A modern proof of this result appear in \cite{BV}.

A morphism  between the crossed modules $\G$ and
$\G'=(G',E',\d',\t')$ is given by a pair $(\f,\p)$ of group morphisms $\f\colon G\to G'$ and
$\p\colon E \to E'$, making the diagram
\begin{equation*}
\begin{CD}
E @>\p>> E' \\
@V\d VV  @VV\d' V\\
 G @>>\f > G'
\end{CD}
\end{equation*}
commutative. 
In addition we must have:  $$\f(X)\t' \p (e)=\p(X \t e); \forall X \in G, \forall e \in E.$$
Crossed modules and their morphisms form a category. This category is a category with colimits, see \cite[3.5]{BV} or \cite{BH3}.

\begin{Example}\label{base}
Let $(M,N,*)$ be a pair of based path connected spaces. If $*'\in N$ is another base point, and $\g$ is a path connecting {$*$ to $*'$,} then there exists a natural isomorphism $\Pi_2(M,N,*) \to \Pi_2(M,N,*')$, where the maps on {homotopy} groups are the usual ones constructed from the path $\g$.
\end{Example}

\begin{Example}\label{freeprod}
Let $\G=(G,E,\d,\t)$ and $\G'=(G',E',\d',\t')$ be crossed modules. The free product  $\G \vee \G'$  of $\G$ and $\G'$  is the pushout, in the category of crossed modules, of the diagram:
\begin{equation*}
\begin{CD}
(\{1\},\{1\},\t,\d) @>>> \G \\
@VVV  \\
\G'
\end{CD}.
\end{equation*}
\end{Example}
Recall that the category of crossed modules is a category with colimits.

\begin{Example}\label{REFER1}
Let $M$ and $N$ be CW-complexes with  unique 0-cells, which we take to be their base points. We have:
$$\Pi_2\big ((M,M^1,*) \vee (N,N^1,*)\big) \cong \Pi_2(M,M^1,*) \vee \Pi_2(N,N^1,*).$$ 
\end{Example}
This follows immediately from the (R. Brown and P.J. Higgins') {2-dimensional} van Kampen Theorem stating that, under mild conditions, the crossed  module functor preserves colimits, and in particular pushouts; see \cite{B2,BV,BH}. This is the case for the wedge product of two well pointed based spaces; see \cite[Theorem 7.1]{BH3}.

\begin{Example}\label{Relations}
Let $\G=(G,E,\d,\t)$ be a crossed module. Suppose that  the elements $a_1,{\ldots} ,a_n\in E$ are such that $\d(a_k)=1_G,k=1,{\ldots} ,n$. Let $F$ be the subgroup of $E$ generated by the elements of the form $X \t a_k$ where $X \in G,k=1,{\ldots} ,n$, thus $F$ is normal in $E$ by the second condition of the definition of crossed modules. {In fact $F$ commutes with all of $E$.} Obviously both $\d$ and the action of $G$ on $E$ descend to $E/F$. Denote the induced map and action  by $\d'$ and $\t'$. It is easy to show that $(G,E/F,\d',\t')$ is a crossed module.
\end{Example}
{We will go back to this construction below.}

\subsubsection{{Presentations of Crossed Modules}}\label{generatorsrelations}
{For details on free crossed modules see \cite{BV,BHu}, for example.}
Let $G$ be a group and $K$ be a set. The free crossed module ${\cal F}\left (\d_0\colon K \to G \right)$  on a map   $\d_0\colon K \to G$ has $G$ {as}  base group. The principal group $E$ is the {quotient of the free group on the set of pairs $(X,m)$, where $X \in G$ and $m \in K$, by the relations:}
\begin{equation}\label{simp}
(X,m)(Y,n)(X,m)^{-1}=(X\d_0(m)X^{-1}Y,n); X,Y \in G;m,n \in K.
\end{equation}
The boundary map $\d \colon  E \to G$ is defined   on generators by $\d(X,m)=X \d_0(m)X^{-1}$, whereas the action of $G$ on $E$ is given simply {by} $X \t (Y,m)=(XY,m)$; here $X,Y \in G$ and $m\in K$. Note that the map {$\i\colon K \to E$} such that $\i(m)=(1_G,m)$ is injective. From now on we consider {$K$} to be included {in} $E$ in this way. 

 The free crossed module on a map $\d_0\colon K \to G$ is defined{, up to isomorphism,} by the following natural  universal property:

\begin{Lemma}
 Let $G'$ be a group and let $\f{\colon} G \to G'$ be a morphism. Let also $\G'=(G',E',\d',\t')$ be a crossed module. Consider a map $\p_0: K  \to E'$ such that $\d'\circ \p_0 =\f \circ \d_0$.
There exists a unique group morphism $\p\colon E \to E'$ extending $\p_0$ in such a way that $(\f,\p)$ is a morphism of crossed modules $\G \to \G'$.
\end{Lemma}

{{Let $G$ be a group and $K$ be a set. Let also $\d_0\colon K \to G$ be a} {map.}  Consider the free crossed module ${\cal F}\left (\d_0\colon  K \to G\right )=(G,E,\d,\t)$.  {A {2-relation}  is, by definition, an element $r$ of $E$ with $\d(r)=1_G$.}}

Let $R=\{r_1,{\ldots} ,r_n\}$ be a set of {2-relations}{, which we can take to be infinite.} {The crossed module presented by the map $\d_0\colon K \to G$, with  2-relations $\{r_1,{\ldots} ,r_n\}$, say ${\cal U}\left (\d_0\colon K \to G;r_1,{\ldots} ,r_n{= 1}\right)$, is,} by definition, the crossed module constructed from  ${\cal F}\left ( \d_0\colon K \to  G\right) $ and $\{r_1,{\ldots} ,r_n\}$ as in {Example} \ref{Relations}.

\subsubsection{A Decomposition of Certain Free Crossed Modules}\label{adec}

This section will only be used for the explicit calculation of the second homotopy  {groups} of the Spun Trefoil {and the Spun Hopf Link} {complements}.

Let $G$ be some group. Let $K$ be a set provided with a map $\d_0\colon  K \to G$. Suppose that $K$ is the disjoint union of the sets $K_1$ and $K_2$. Let $\d^1_0$ and $\d^2_0$ be the restrictions of $\d_0$ to $K_1$ and $K_2${, respectively}. 
Suppose also that $\d_0^2(k)=1_G, \forall k \in K_2$. Let $F_2$ be the subgroup of the principal group  of ${\cal F}\left ( \d_0:K \to G \right)=(G,E,\d,\t)$ generated  by the elements $X \t f$, where $f \in K_2$ and $X \in G$. Then it is easy to see that $F_2$ commutes with $E$ and that it is the free {abelian} module {over} ${\rm coker}(\d)$ with base $K_2$, with the obvious action of $G$. This follows from relations  (\ref{simp}) and the fact that for any crossed module $(G',E', \d',\t')$, the image  ${\rm im} (\d')$ of $\d'$ acts trivially on ${\rm ker}(\d')$; {a consequence of} the second condition of the definition of crossed modules. 

In fact, the group $F_2$ has an algebraic complementary $F_1$ in $E$, as a group and as a $G$-module,  where $F_1$ is isomorphic with the principal group of the free crossed module on the map $\d_0^1\colon K_1 \to G$. This follows immediately from  relations  (\ref{simp}). 

Therefore we have:
\begin{Proposition}
Let $G$ be a  group, and  let $K$ be a set provided with a map $\d_0\colon  K \to G$. Suppose that $K$ is the disjoint union of the sets $K_1$ and $K_2$. Let $\d^1_0$ and $\d^2_0$ be the restrictions of $\d_0$ to $K_1$ and $K_2${, respectively}. 
Suppose also that $\d_0^2(k)=1_G, \forall k \in K_2$.  Then as a $G$-module, and as a group, the principal group of the free crossed module {$(G,E,\d,\t)$} on $\d_0\colon  K \to  G$ is the direct sum of the principal group {$F_1$} of the free crossed module on  ${{\d_0^1}\colon  K_1 \to G}$  and the free ${\rm coker} (\d)$-module {$F_2$} with base {$K_2$}{, with {the} obvious action of $G$. This direct sum has the natural boundary map to $G$ where $\d(F_2)=\{1_G\}$}. 
\end{Proposition}
{This result is also valid, with the obvious modifications, for the case of a  crossed module ${\cal U} \left (\d_0 \colon K \to G;r_1,...,r_n=1\right)$ presented by a map $\d_0\colon K \to G$ with  2-relations $r_1,...,r_n$, as long as each 2-relation is contained in either  $F_1$ or $F_2$.}

\subsection{{The Significance of the Fundamental Crossed Module $\Pi_2(M,M^1)$}}

Whenever  $E$ is an abelian group and we have a left action of  $G$  on $E$ by automorphisms, then $(G,E,\d=1_G, \t)$ is always a crossed module.
For any based path connected topological space $M$,  the group $\pi_2(M)$ is abelian and $\pi_1(M)$ acts on $\pi_2(M)$ by automorphisms. Therefore we have a crossed module $\pi_{1,2}(M)$ for any based topological space. See for example  \cite{L,M} for calculations of $\pi_{1,2}(M)$ when $M$ is the complement of a knotted surface in $S^4$.

A first idea {about} how to employ the notion of a  crossed module to define
invariants of manifolds could be to consider the crossed module $\pi_{1,2}(M)$.  However, even
when $M$ is a compact manifold, it is not  certain that $\pi_2(M)$ is
finitely generated as a module over $\pi_1(M)$; see  \cite[problem 5]{L}  and, less directly related, {\cite[problems 6 and 13]{L}} for the important case of complements of knotted surfaces. Therefore,  $\pi_{1,2}(M)$  is not a very practical invariant since it not easy to distinguish {between} two non-finitely  generated $\pi_1(M)$-modules.

Another solution is to consider the more tractable relative case. Let $M$ be a CW-complex, and let $M^1$ be its 1-skeleton. Consider the  crossed module  
$\Pi_2(M,M^1)=\left (\pi_1(M^1),\pi_2(M,M^1),\d,\t\right)$.
  Despite the asymmetry introduced by choosing a particular 1-skeleton of $M$, this crossed module determines $\pi_2(M)$ and $\pi_1(M)$,  which fit inside the exact sequence:
\begin{equation}\label{kernel}
\{0\} \to \pi_2(M) \to \pi_2(M,M^1) \ra{\d} \pi_1(M^1) \to  \pi_1(M)\to \{1\},
\end{equation}
since $\pi_2(M^1)=\{0\}$.  In fact the crossed module $\Pi_2(M,M^1)$  determines also  the  $k$-invariant $k(M) \in H^3(\pi_1(M),\pi_2(M))$. The  group cohomology class $k(M)$ is determined from the classical correspondence between 3-dimensional group cohomology classes and crossed modules; see for example \cite{B2,ML,KB,JH}. Therefore  the crossed module $\Pi_2(M,M^1)$ determines the topological 2-type of $M$, see \cite{MLW}; {and thus it is} strictly stronger than $\pi_1(M)$ and $\pi_2(M)$ alone; see \cite{Bau,Lo}.

Choosing the apparently less charming relative case is also justified by the fact that $\Pi_2(M,M^1)$ does not depend on the cellular decomposition of $M$ up to free products{, in the category of crossed modules,} with $\Pi_2(D^2,S^1)$, as we will prove below in \ref{mainset}. This expands an old result of {J.H.C.} Whitehead stating that the homotopy type of  the crossed module $\Pi_2(M,M^1)$ depends only on the homotopy type of $M$, as a space; see \cite{W3,BH2}.

Therefore the fundamental crossed module  $\Pi_2(M,M^1)${, up to free products with $\Pi_2(D^2,S^1)$,} is a genuine homotopy invariant.  This {allows} us to obtain a homotopy invariant of finite connected CW-complexes for any finite crossed module, see \ref{CMI}.

 We will see {below} that in {the} case when $M$ is a  CW-complex and $M^1$ is its 1-skeleton then $\Pi_2(M,M^1)$ is, in principle,  calculable. 

\subsubsection{ Whitehead's Isomorphisms}\label{wi}
For a more complete treatment of these issues see \cite{Br}.
Let $M$ be a path-connected topological space with a base point $*$. Let $N$ be a topological space obtained from $M$ by attaching some 2-cells (or 2-handles) $s_1,{\ldots} ,s_n$. Choose a base point {$*_i$  on the boundary of  each  2-cell, where $i=1,\ldots,n$}.  If we are provided paths connecting $*_i$ with $*$, for $i=1,{\ldots} ,n$,  then each 2-cell $s_i$ can be identified  uniquely with  an element of  $\pi_2(N,M,*)$. This does not depend on the path chosen up to acting  by some element of {$\pi_1(M,*)$.}

 Consider the map {$\d_0\colon \{s_1,{\ldots} ,s_n\}\to {\pi_1(M,*)}$} induced by the {attaching maps} of each cell $s_i$. This map (well defined up to {conjugations by elements} of ${\pi_1(M,*)}$) is also fixed by the chosen paths  connecting $*_i$ with $*$, for $i=1,{\ldots} ,n$.

Recall the following theorem, due to J.H.C. Whitehead. For the original
proof see \cite{W1,W2,W3}; see also \cite[5.4]{BV} and \cite{B3,BH2,GH}.
 
\begin{Theorem}  \label{Whitehead} 
The natural morphism from the free crossed module on the  map {$\d_0\colon \{s_1,\ldots,s_n\}\to \pi_1(M,*)$} into $\Pi_2(N,M,*)$ is an isomorphism of crossed modules.
\end{Theorem}
This result (usually called Whitehead's Theorem) is one of the most important results that we will use.

Let $M$ be a CW-complex.  For each $k \in \N$, let $M^k$ denote the $k$-skeleton of $M$. For simplicity, suppose that $M$ has a unique $0$-cell, which we take to be  its base point $*$. 
Note that from  the Cellular Approximation Theorem we have $\Pi_2(M,M^1,*) \cong \Pi_2(M^3,M^1,*)$. 

Consider the group complex: 
$${\ldots}  \ra{\d_4} \pi_3 (M^3,M^2,*) \ra{\d_3} \pi_2(M^2,M^1,*) \ra{\d_2} \pi_1(M^1,*)\ra{p} \pi_1(M,*),$$
with the obvious boundary maps. This is a crossed complex of free type, called the {``Fundamental Crossed Complex of $M$''}; see for example \cite{B2,B3}. In particular, we have an action of the group $\pi_1(M^1,*)$ on all the {other groups}, preserving the {boundary maps}, and such that, if $n>2$, then the action of $\pi_1(M^1,*)$ on $\pi_n(M^n,M^{n-1},*)$ factors through {the projection map} $p\colon \pi_1(M^1,*)\to \pi_1(M,*)$.

 {For each $n$-cell $c^n$ of $M$, where $n>1$, choose a base point {on the boundary of it,} as well as a path from the cell base point to $*$. Therefore $c^n$ determines an element  $c^n \in\pi_n(M^{n},M^{n-1},*)$, and its {attaching} map can also  be identified with an {element of $\pi_{n-1}(M^{n-1},*)$. If $n>2$ the projection of this element in $\pi_{n-1}(M_{n-1},M_{n-2},*)$ is the boundary map $\d_n(c^n)\in \pi_{n-1}(M_{n-1},M_{n-2},*)$ of $c^n$, considering the fundamental complex complex of $M$.}

It is well known, see \cite[V.1]{GW}, \cite{W3} or \cite{B2} that{, if $n>2$, then} the natural map from the free $\Z(\pi_1(M,*))$-module over the $n$-cells of $M$ into $\pi_n(M^n,M^{n-1},*)$ is an isomorphism. This is also a result of J.H.C. Whitehead. 
In particular, the group $\pi_3 (M^3,M^2,*)$ is the free $\Z(\pi_1(M^2,*))$-module on the group elements defined by the $3$-cells of $M$. 

 From the homotopy exact sequence of the triple $(M^3,M^2,M^1)$ it follows that:
\begin{Lemma}
$$\pi_2(M,M^1,*)=\pi_2(M^2,M^1,*)/{\rm im}(\d_3).$$
\end{Lemma}
Let $c_1^3,{\ldots} ,c_{n_3}^3$ be the 3-cells of $M$. Each one of them defines an element of $\pi_3(M^3,M^2,*)$, as well as its boundary ${\d_3(c^3_k)}\in\pi_2(M^2,M^1,*),k=1,{\ldots} ,n_3$, well defined if we make the choices above. Let also ${\{c_1^2,{\ldots} ,c_{n_2}^2\}}$ be the set of 2-cells {of $M$}, where $c_k^2$ attaches along $\d_2(c_k^2) \in \pi_1(M^1,*),{k=1},{\ldots} ,n_2$.

\begin{Theorem}\label{Attach3}{
The crossed module $\Pi_2(M,M^1,*)$ is the  crossed module presented by the  map  from the  set of 2-cells {of $M$} into {$\pi_1(M^1,*)$,} defined from the attaching maps of each 2-cell of $M$, with one {2-relation} for each 3-cell of $M$. More precisely:}
$${\Pi_2(M,M^1,*)\cong {\cal U} \left ( {\{c_1^2,{\ldots} ,c_{n_2}^2\}}\ra{\d_2} {\pi_1(M^1,*)}; \d_3(c^3_1),{\ldots} ,\d_3(c^3_{n_3}){=1}\right ).}$$
Note that $\d_2 \circ \d_3=1$.
\end{Theorem}

\begin{Proof}
By Whitehead's Theorem and the previous lemma,  we only need to prove that ${\rm im} (\d_3)$ is the subgroup of $\pi_2(M^2,M^1,*)$ generated by the elements $X \t \d_3(c^3_i)$, where $X \in \pi_1(M^1,*)$ and $i=1,{\ldots} ,n_3$. This follows from the fact that $\pi_3(M^3,M^2,*)$ is the free $\Z\big (\pi_1(M^2,*)\big )$-module on the set of 3-cells of $M$.
\end{Proof}

This theorem tells us  that, in principle, if $M$ is a CW-complex and $*$ is a 0-cell {of $M$}, then the  crossed module $\Pi_2(M,M^1,*)$ can be calculated. The only possible difficulty is the determination of the boundary maps in the fundamental crossed complex of $M$. This can be solved for example for simplicial complexes by the homotopy addition lemma in \cite[page 175]{GW}. The case of {complements of} knotted surfaces can also be solved by using a particular handle decomposition of them. This is the main aim of this article. 

\subsubsection{{The Dependence of $\Pi_2(M,M^1)$ on the Cell Decomposition of  the CW-complex $M$}}\label{mainset}
Let $(N,M)$  be a pair of {connected} CW-complexes such that the inclusion of $M$ in $N$ is a homotopy equivalence. Let $M^1$ and $N^1$ be, respectively, the 1-skeletons of $M$ and $N$. Suppose that $M$ has a unique 0-cell, which we take to be the base point $*$ of $M$ and $N$, so that both $M$ and $N$ are well pointed.

The group $\pi_1(M^1,*)$ is the free group on the set $\{d_1,{\ldots} ,d_m\}$ of 1-cells of $M$. There exist also $c_1,{\ldots} ,c_n \in \pi_1(N^1,*)$ such that $\pi_1(N^1,*)$ is the free group $F(d_1,{\ldots} ,d_m,c_1,{\ldots} ,c_n)$ on the set    $\{d_1,{\ldots} ,d_m,c_1,{\ldots} ,c_n\}$. {These elements of $\pi_1(N^1,*)$ define elements of $\pi_1(N,*)$ in the obvious way.}

\begin{Theorem}
There exists a homotopy equivalence:
$$(N,N^1,*)\cong (M,M^1,*) \vee (D^2,S^1,*)^{\vee n}.$$                   
 \end{Theorem}

\begin{Proof} \footnote {This argument arose in a discussion with Gustavo Granja.}  
Since $M$ is a subcomplex of $N$, and $N$ is homotopic to $M$, it follows that $M$ is a strong deformation retract of $N$. By the Cellular Approximation Theorem, we can suppose, further, that there exists a retraction $r\colon N \to M$ sending $N^1$ to $ M^1$, and such that $r \cong \id_N$, relative to $M$.  In particular if $k \in \{1,{\ldots} ,n\}$ then we have {that} $c_k r_*(c_k)^{-1}=1_{{\pi_1 (N,*)}}$, {in $\pi_1(N,*)$.}
Define a map $$f\colon (P,P^1,*) \doteq (M,M^1,*) \vee \bigvee_{k=1}^n (D^2_k,S^1_k,*)\to (N,N^1,*)$$ in the following way.

 First of all, send $(M,M^1,*)$ identically to its copy $(M,M^1,*)\subset (N,N^1,*)$. Then we send each $S^1_k$ to the element $c_k r_*(c_k)^{-1}\in \pi_1(N^1,*)$. Since these elements are null homotopic in $(N,*)$, this map extends to the remaining 2-cells of $(P,P^1,*)$.

 Let us prove that $f\colon (P,P^1,*) \to (N,N^1,*)$ is a homotopy equivalence. It suffices to prove that $f\colon (P,*) \to (N,*)$ and ${f^1 \doteq} f_{|P^1}\colon (P^1,*) \to (N^1.*)$ are based homotopy equivalences, since all inclusion maps are cofibrations{; see}  for example {\cite[6.5]{M2}.} Notice that the  result proved there is also valid in the base case, as long as all the spaces considered are well pointed, which is the case {here}.

{It is immediate that $f$ is an equivalence of homotopy $(P,*) \to (N,*)$, since $f$ extends the inclusion map $(M,*) \to (N,*)$.} Let us show that ${f^1\colon (P^1,*) \to (N^1,*)}$ is a  homotopy equivalence. It is enough to prove that the induced map ${f^1_*}\colon  \pi_1(P^1,*) \to \pi_1(N^1,*)$ is an isomorphism. Note that  $\pi_1(P^1,*)$ is (similarly with $\pi_1(N^1,*)$) isomorphic  with the free group on the set $\{d_1,{\ldots} ,d_m,c_1,{\ldots} ,c_n\}$. The induced map on the fundamental groups has the form: 
$${f^1_*}(d_k)=d_k, k=1,{\ldots} ,m, \textrm{ and }{f^1_*}(c_k)=c_kr_*(c_k)^{-1}, k=1,{\ldots} ,n.$$
 Notice that $r_*(c_k) \in F(d_1,{\ldots} ,d_m), k=1,{\ldots} ,n$. Consider the morphism $g$ of $F( d_1,{\ldots} ,d_m,c_1,{\ldots} ,c_n)$ on itself such that:
$$g(d_k)=d_k, k=1,{\ldots} ,m \textrm{ and } g(c_k)=c_kr_*(c_k), k=1,{\ldots} ,n.$$
 Therefore $({f^1_*} \circ g )(d_k)=d_k, k=1,{\ldots} ,m$ and 
\begin{align*}({f^1_*} \circ g)(c_k)&={f^1_*}(c_k r_*(c_k))\\&={f^1_*}(c_k) {f^1_*}(r_*(c_k))\\&=c_kr_*(c_k)^{-1}r_*(c_k)\\&=c_k,k=1,{\ldots} ,n.\end{align*}
 On the other hand $(g \circ {f^1_*})(d_k)=d_k, k=1,{\ldots} ,m$ and  
\begin{align*}
(g \circ {f^1_*})(c_k)&=g(c_k r_*(c_k)^{-1})\\
&=c_k r_*(c_k) g(r_*(c_k^{-1}))\\
&= c_k r_*(c_k) r_*(c_k^{-1})\\
&=c_k,k=1,{\ldots} ,n.
\end{align*}
 This  proves that $g^{-1}={f^1_*}$, which finishes the proof.
\end{Proof}

\begin{Corollary}\label{referee}
Let $M$ and $N$ be  CW-complexes with  unique $0$-cells, which we take to be their base points $*$ and $*'$. Suppose that $M$ and $N$ are homotopic as spaces. There exists $m,n \in \N$ such that:

$$(M,M^1,*)\vee (D^2,S^1,*)^{\vee n} \cong (N,N^1,*')\vee (D^2,S^1,*')^{\vee m},$$
thus, in particular:
$$\Pi_2(M,M^1,*) \vee \Pi_2(D^2,S^1,*)^{\vee n} \cong \Pi_2(N,N^1,*')\vee \Pi_2(D^2,S^1,*')^{\vee m},$$
where the free product is taken in the category of crossed modules.
\end{Corollary}
The second part of this result was suggested by a referee of a previous version of this article.

\begin{Proof}
The {pointed} spaces $(M,*)$ and $(N,*')$ are homotopic. Let $P$ be the reduced mapping cylinder of some {pointed} homotopy equivalence $(M,*) \to (N,*')$, chosen to be cellular. Therefore both $M$ and $N$ are cellularly included in $P$ (provided with its usual CW-decomposition), and they intersect along their {base} points, both coinciding with the unique 0-cell of $P$. 
The first result follows from the previous theorem and the fact that the inclusions of $M$ and $N$ {in} $P$  are homotopy equivalences. The second one follows directly from the 2-Dimensional van Kampen Theorem; see \cite{B2,BH,BV}. 
\end{Proof}

Note that we necessarily have $b_1(M^1)+n=b_1(N^1)+m=b_1(P^1)$. We will need to use this fact later. {Here, if $K$ is a CW-complex, then $b_1(K)$ denotes the first Betti number of it. In addition $P$ is the reduced mapping cylinder of some {pointed} homotopy equivalence $(M,*) \to (N,*')$, chosen to be cellular.}

It is easy to see that $\Pi_2(D^2,S^1)=(\Z,\Z,\id,\t)$, where $\t$ is the trivial action. We can prove {this} from the  {long homotopy} exact sequence of $(D^2,S^1)$, or alternatively by using Whitehead's Theorem together with the explicit description of free crossed modules in \ref{generatorsrelations}.

{The results of this subsection extend in a natural way for the case of the fundamental crossed complex of a CW-complex; see \cite{FM2}.}

\section{Complements of Knotted Surfaces}

For the case of complements of knotted surfaces {in $S^4$}, it is {convenient to work with handle decompositions,}  more flexible than CW-decompositions, considered in the previous section.

 Let $M$ be a manifold with a handle decomposition,  and let $M^{(n)}$ be the handlebody {made out of} the handles of $M$ of index smaller or equal to $n$ (the $n$-handlebody of $M$).  It is well known that a handle decomposition of  $M$  determines a topological space  $\hat{M}$ of the same homotopy type of $M$, with a CW-decomposition where each $n$-handle of the manifold $M$ generates an $n$-cell of the CW-complex $\hat{M}$, see \cite[Chapter 6]{RS} or  \cite [Chapter III]{Mz}, for example. Intuitively, $\hat{M}$  is obtained from $M$ by shrinking any $n$-handle to  an $n$-cell going along its core.  More precisely, it is possible to prove that{, if $M$ is a manifold with a handle decomposition,} then there exists a homotopy equivalence $M \to \hat{M}$ preserving the filtrations of $M$ and $\hat{M}$ given by their handle and cell decompositions, respectively, and such that the restriction  maps $M^{(n)}\to {\hat{M}^n}$ are  homotopy equivalences for each $n$; see \cite[Proposition 3.4]{Mz}. Given that  all inclusion maps are cofibrations, it thus follows that there exists  a filtered homotopy equivalence $M \to \hat{M}$.

 Let again $M$  be a manifold with a handle decomposition.  Suppose that the base point of $M$ is contained in one of the 0-handles of $M$. The above equivalence provides a homotopy equivalence $(M,M^{(1)},*) \cong (\hat{M},\hat{M}^1,*)$. Here the base point of $\hat{M}$ is taken to be one of the 0-cells of $\hat{M}$. 
 Therefore, the  results of the previous chapter apply with the obvious adaptations   to the crossed module $\Pi_2(M,M^{(1)})$.

\subsection{ {Handle Decompositions} of  Complements of Knotted Surfaces} \label{handledec} 
 For details on knotted surfaces, in particular movie presentation of them, we refer the reader to \cite{CS,CRS,CKS,F}. We work in the smooth category.

\subsubsection{Movies of Knotted Surfaces}

Let $\S\subset S^4=D^4_- \cup ( S^3 \times I) \cup D^4_+$,  where from now on $I=[-2,2]$, be a knotted surface. In other words $\S$ is  a (locally flat) embedding of a closed  2-manifold {into $S^4$.} We want to calculate $\Pi_2(M,M^{(1)})$, where $M = S^4\setminus \n(\S)$. Here $\n(\S)$ is an (open) regular neighborhood of $\S$ in $S^4$ and   $M^{(1)}$ is the handlebody made {out of} the 0- and 1-handles of a handle decomposition of $M$ (the 1-handlebody {of} $M$).   This will provide a description of the algebraic 2-type of the complement of ${\n(\S)}${, in the form of a crossed module.}

We need to  construct  a  handle decomposition of the complement of $\n(\S)$ in $S^4$. Such handle decomposition can be defined from a movie of $\S$; {see \cite{G,GS,CKS}.}

 Up to isotopy, we can suppose that $\S \subset S^3 \times I$. Suppose that  the projection on  $I$ defines a Morse function on the knotted surface $\S$. In particular, for each non critical $t\in I$, the set  $\S_t =\S \cap \left (S^3 \times \{t\}\right)$ is a link in $S^3$ (a still of $\S$). Between critical values, the link $\S_t$ will undergo an isotopy of $S^3$. At critical points of index $0$, $1$ or $2$, the link $\S_t$ will {go through}  Morse modifications, called, respectively, ``Minimal Points'', ``Saddle Points'' and ``Maximal Points''; see figure \ref{crit1}. The 1-parameter family of links $t \mapsto \S_t$, with the modifications at non-generic points,  will define what is  called a ``movie'' of the knotted surface $\S \subset S^4$.

\begin{figure}
\begin{center}
\includegraphics{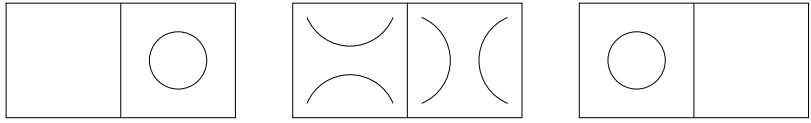}
\end{center}
\caption{Morse Modifications, respectively a ``Minimal Point'',  a ``Saddle Point'' and  a ``Maximal Point''. Each of these modifications goes from left to right.}
\label{crit1}
\end{figure}

\begin{figure}
\begin{center}
\includegraphics{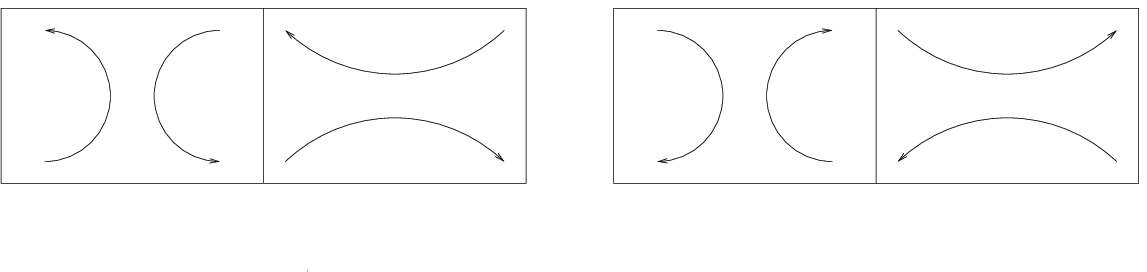}
\end{center} 
\caption{Oriented saddle point transitions.}
\label{osaddle}
\end{figure}
 If the knotted surface $\S$ is oriented, then each link $\S_t$ appearing in the movie (now called an oriented movie) of $\S$ will have a natural orientation.  See figure \ref{osaddle} for the oriented version of the  saddle point moves. Any oriented movie defines an oriented knotted surface, up to isotopy.

\subsubsection{Hyperbolic Splittings and Knots with Bands}
Let $\S\subset S^4=D^4_- \cup ( S^3 \times I) \cup D^4_+$  be a knotted surface presented by a movie $t \mapsto \S_t$. {Recall that $I=[-2,2]$.}  By using isotopy, we can suppose that $\S \subset S^3 \times [-1,1]$, and, moreover, that all minimal points occur in $S^3 \times \{-1\}$, all maximal points occur in $S^3 \times \{1\}$, and all saddle points occur in $S^3\times \{0\}$; see for example \cite[Chapter 1]{CKS} or \cite{KSS}. This is a well know result. This type of movies of knotted surfaces are usually called ``hyperbolic splittings''.

Consider a knotted surface $\S$ represented by a hyperbolic splitting.  Therefore, for all $t$ in $(-1,0)$,  the still $\S_t=\S \cap \left( S^3 \times \{t\}\right)$ of $\S$  will be an unlink with a fixed number of  components, and $\S_t$ will undergo an isotopy of $S^3$ in this interval. The same is true for $t \in (0,1)$. At $t=0$ the link $\S_t$ will undergo saddle point transitions; see figure \ref{osaddle}. {We have in addition minimal and maximal points at $t=-1$ and $t=1$, respectively.}

All this information used for  constructing  a knotted surface is highly redundant, \cite{L}. In fact,  the  knotted surface constructed in this way depends only on  the saddle point transitions at $t=0$, as well as the configuration immediately before and after $t=0$; see for example \cite{KSS}  for a proof.

The saddle point transitions which happen  at $t=0$ can be encoded by a knot with bands; see \cite{KSS,SW,CKS}.  Another usual (and equivalent) presentation is to use marked vertex diagrams; see \cite{L,CKS,Yo}. The former are more useful for our purposes since with them we can suppose that the configuration  immediately after the saddle points is a  standard diagram of the unlink.

 A knot with bands  is a knot together with some $I\times I$ bits, intersecting
 the knot along $\d I \times I$. A knot with bands is said to be oriented if we have orientations on the thin edges of it, having the configuration of  figure \ref{kwb}, or its mirror image, at the edges incident to a band. Henceforth, all knots with bands will be oriented.  A knot with bands $K$ determines two oriented knots $K_+$ and $K_-$ called the post-knot and pre-knot of $K$, see figure \ref{prepos}. Note that our convention is opposite to the one in \cite{SW}.
\begin{figure} 
\begin{center}  
\includegraphics{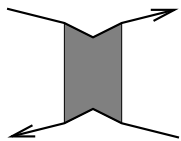} 
\end{center} 
\caption{A bit of a knot with bands.}
\label{kwb} 
\end{figure} 
\begin{figure}
\begin{center} 
\includegraphics{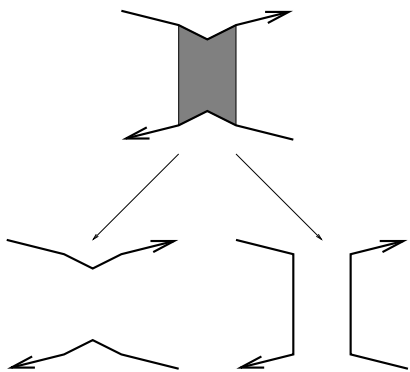} 
\end{center}
\caption{A piece of a knot with bands together with its post-knot and pre-knot (from
  left to right).}
\label{prepos}
\end{figure}

Let $K$ be a knot with bands such that the post- and {pre-knots} of $K$ are oriented unlinks.  We can construct an oriented  knotted surface by choosing an isotopy from the post-knot $K_+$ of $K$ to the standard unlink diagram for $t>0$, and by capping all the circles of it in the obvious way, and analogously for the interval $[-1,0)$. The final result does not depend on the choices made, up to isotopy. {In fact, it depends only on the isotopy class of $K$.}

In this article we will consider {this} description of knotted surfaces. {However, we will  need to use the movie picture associated  a presentation of a knotted surface by  knot with bands when constructing the handle decomposition of the complement.}

For a  set of moves relating any two knot with bands {representations} of the same knotted surface (up to isotopy), we refer the reader to \cite{SW}. We will not need to use that result.

\subsubsection{Construction of the Handle Decomposition}

Let $\S\subset S^4=D^4_- \cup ( S^3 \times I) \cup D^4_+$ be a knotted surface presented by a movie $t\mapsto \S_t$. Therefore $\S$ is provided {with} a Morse function (the projection on $I$) and, away from critical points, $\S_ t$ is a link in $S^3$. There exists a natural  handle decomposition of the complement $M$ of an open  regular neighborhood $\n(\S)$ of $\S$ (in principle defined up to handle-slides and isotopy)  where minimal/maximal points will induce $1/3$-handles of the decomposition, and saddle points induce $2$-handles; see \cite[section 6.2]{GS}, \cite[3.1.1]{CKS}, or \cite{G}. This is very easy to visualize in dimension $3$. To calculate the fundamental crossed module $\Pi_2(M,M^{(1)})$, however, we need an explicit description of this handle decomposition. We follow now \cite[3.1.1]{CKS} and \cite[6.2]{GS}, where the  missing bits of our description can be found.

Let $K$ be a knot with bands, representing the knotted surface $\S${, thus $\S \subset S^3 \times [0,1]$}. Choose a regular projection of $K$. Let $t \mapsto \S_t$ be an associated movie for $\S$, a hyperbolic splitting.  For simplicity (and without loss of generality), we will suppose that the post-knot of $K$ is a standard  diagram of the unlink (a disjoint union of unknotted circles). This will fix a handle decomposition of the complement {$M$} of a regular neighborhood of $\S$, up to isotopy.  From now on  we will suppose that all knots with bands representing knotted surfaces are on this form.

 For any set $A$ in $ S^4=D^4_- \cup ( S^3 \times I) \cup D^4_+$, let $A_t=A \cap (S^3 \times \{t\})$ and $A_{\leq t}=A \cap \big ((S^3 \times [-2, t]) \cup D^4_-\big )$. Then there exists some small  $\e >0$ such that the topology of $M_{\leq t}$ does not change in the intervals $[-2, -1-\e]$, $[-1+\e,-\e]$, $[\e,1-\e]$ and  $[1+\e, 2]$. In between these intervals, the manifold 
$M_{\leq t}$ will undergo attachment of handles of index, respectively, 1,2,3 and 4.

 The manifold $M_{\leq (-1+\e)} \cong M_{\leq -\e}$ is a 4-manifold obtained {from} attaching $n$ $1$-handles to $D^4 \cong M_{\leq (-1-\e)}$. Here $n$ is the number of components of the pre-knot of $K$. In fact a Kirby diagram for  $M_{\leq -\e}$  is obtained from the pre-knot $K_-$ of $K$ by turning the circles of it  (considered to be 0-framed) into dotted circles, in the notation of Kirby, \cite{K}; see \cite[Section 6.2]{GS}.  This is also clear from the construction in \cite[3.1.1]{CKS}.

The manifold $M_{\leq \e}$ is  obtained from $M_{\leq -\e}$ {from} attaching 2-handles, each of which {is} determined by a band of $K$.
For each $t\in [\e,1-\e]$, the manifold $M_{\leq t}$   is therefore a 4-dimensional handlebody made from $0$-, $1$- and $2$-handles, {the 2-handlebody $M^{(2)}$} of the knotted surface complement $M$.  

We can easily obtain a Kirby diagram for $M^{(2)}$  from the knot with bands $K$; cf. \cite[Section 6.2]{GS}. Consider a dotted circle for each circle of $K_-$, the pre-knot of $K$. Then  each 2-handle should be attached along a framed circle $S^1 \times I$ encircling {one of the  bands of $K$,} with  framing parallel to the core of it; see figure \ref{bandkirby}. In particular  each 2-handle attaches along a 0-framed {circle}. To understand the subsequent attachment of 3-handles, it is convenient to draw the framed {circle} $S^1 \times I$ encircling each  band  of $K$ in a way such  that the framing of it goes along     {almost} the entire length of the band, as in figure \ref{bandkirby}.
\begin{figure}
\centerline{\relabelbox 
\epsfysize 8cm
\epsfbox{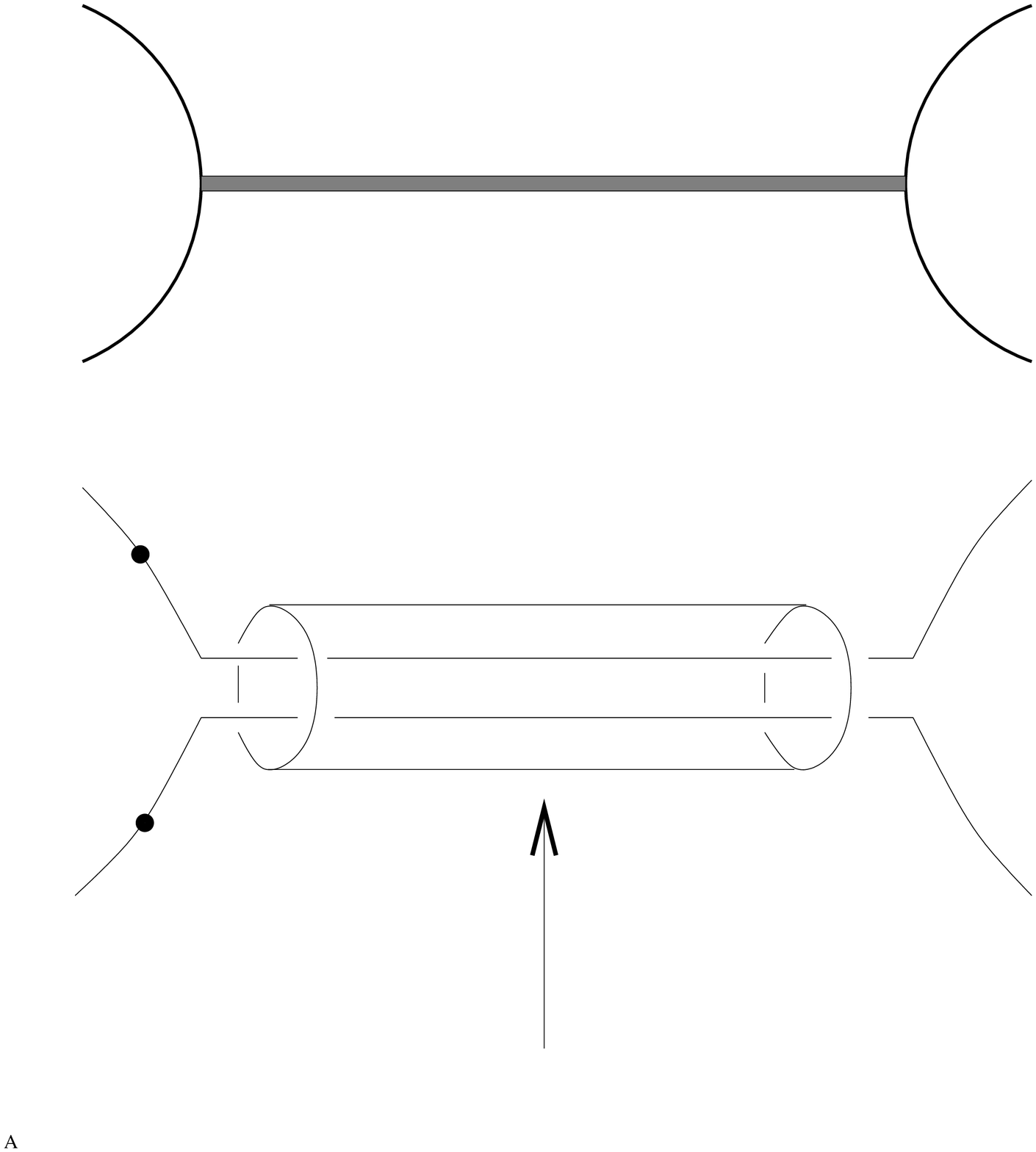}
\relabel{A}{$\s{\textrm{framed circle } S^1 \times I \textrm{ along which a 2-handle attaches}}$}
\endrelabelbox }
\caption{The Kirby diagram for the 2-handlebody  $M^{(2)}$ of the complement of a knotted surface  obtained from a knot with bands. {The original diagram appears on top.} }
\label{bandkirby}
\end{figure} 
 
 To describe the attaching map of each 2-handle (and not simply the attaching  region) we need an orientation of its attaching sphere. Such orientation can be fixed by  an orientation of the core of the associated band; see    figure \ref{bandshandle}. In subsequent  drawings of knots with bands, there will be arrows denoting the orientation of both the thin and fat strands of it; see figure \ref{orientation}. 
 
\begin{figure}
\centerline{\relabelbox 
\epsfysize 5cm
\epsfbox{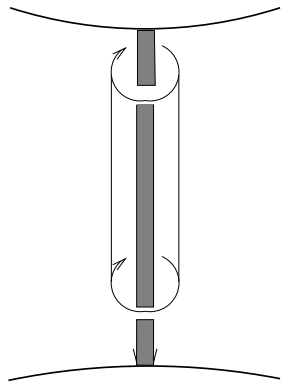}
\endrelabelbox }
\caption{The oriented framed {circle $S^1 \times I$} determined by a band of a knot with bands with its core oriented. }
\label{bandshandle}
\end{figure}

\begin{figure}
\centerline{\relabelbox 
\epsfysize 2cm
\epsfbox{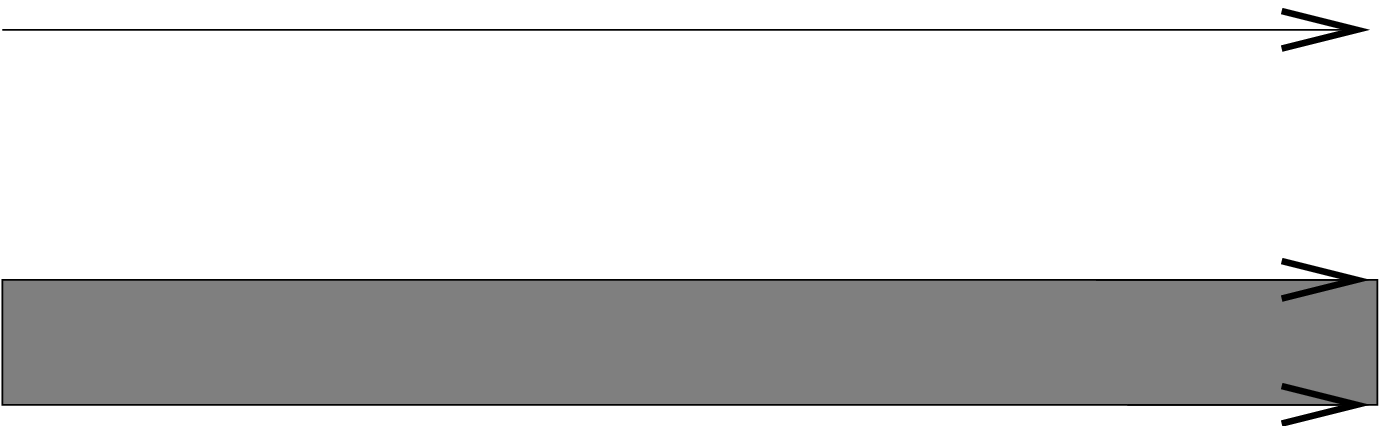}
\endrelabelbox }
\caption{Denoting orientations {of} thin and fat components of a knot with bands. }
\label{orientation}
\end{figure}

The  handles of index 3 attach along regions diffeomorphic with $S^2 \times I$. In the case of complements of knotted surfaces, $3$-handles correspond to maximal points.  The attaching sphere of each 3-handle will be a  sphere $S^2$ containing one of the circles of the post-knot of $K$ {{in the region inside}  it}, as in figure \ref{3handle2}; see \cite[3.1.1]{CKS} and \cite[6.2]{GS}. Recall that we suppose that  the post-knot $K_+$ of $K$ is a standard diagram for the unlink.

 In the Kirby diagram for the knotted surface complement {$M$,} the configuration can be more complicated due to the previous attachment of $1$- and $2$-handles. 
 We can suppose that each of the  attaching spheres $S^2$ determined by the circles of $K_+$ intersects the framed circles $S^1\times I$ determined by the bands of $K$ transversally, thus along a disjoint union of circles $S^1$, {circles} which we can also suppose go around the corresponding band of $K$. This type of intersections will be called essential. Therefore, in the vicinity  of the circles of $K_+$, the configuration of the Kirby diagram {of $M$} will  look like  figure \ref{3handle}.

 Note that since we perform surgery on the framed knots $S^1 \times I$ appearing in figure \ref{3handle}, the shown embedded sphere $S=S^2$ is well defined. In fact the intersection of $S$ with the {previously} attached $2$-handles is, in this particular case,  a disjoint union of three disks  $U_1$, $U_2$ and $U_3${, whose} boundary is the intersection of $S$ with the framed circles $S^1\times I$ determined by the bands; see figure \ref{3handle}. These disks $U_1$, $U_2$ and $U_3$ are parallel to the core of the corresponding 2-handles {of $M$}. This remark  continues to hold for more complicated configurations, with the obvious adaptations. Therefore we have:

\begin{Lemma}\label{B}
The attaching sphere  of each $3$-handle is a sphere $S^2$ containing one, and only one, of the circles of the post-knot $K_+$ of $K$ {{in the region inside}  it}, as in figure \ref{3handle2}. We can suppose that each attaching sphere $S^2$ intersects the framed circles $S^1 \times I$  determining the attachment of $2$-handles transversally, thus each connected component of the  intersection  is a circle $S^1$, which {furthermore} we can suppose is {linking} the associated band {of $K$}, locally (an essential intersection). {Moreover}, the attaching sphere $S^2$ of each 3-handle intersected with the $2$-handles is a disjoint union of  disks (one for each essential intersection), each of which is parallel to the  core of the corresponding 2-handle. {The  boundary of each of these disks is the corresponding connected component of  the intersection of the attaching sphere $S^2$ with the framed circles determining the attachment of 2-handles.} 
\end{Lemma}
Note that the framed circle determined by a band may intersect the attaching sphere for a 3-handle more that once.

\begin{figure}
\centerline{\relabelbox 
\epsfysize 5cm
\epsfbox{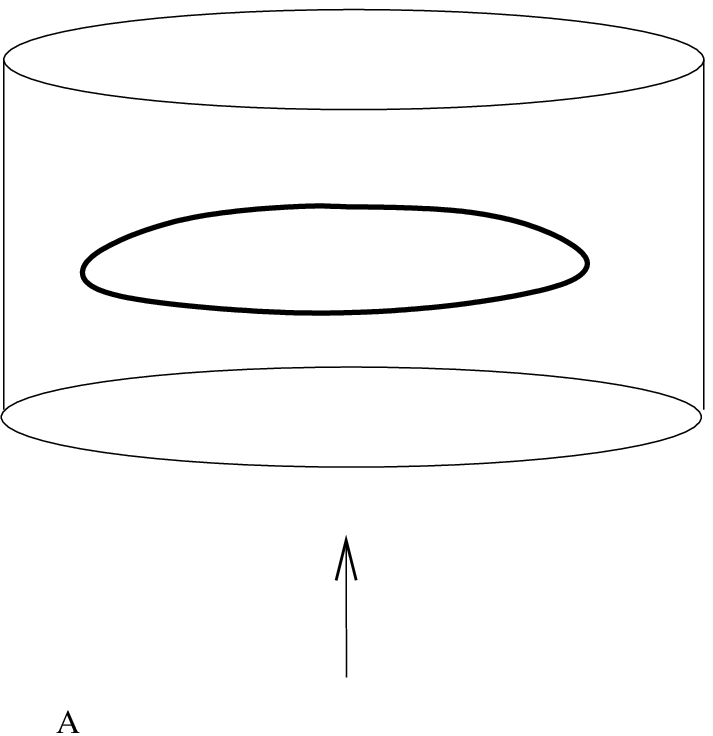}
\relabel{A}{$\s{\textrm{attaching sphere } S^2 \textrm{ for a 3-handle}}$}
\endrelabelbox }
\caption{The attaching sphere for a 3-handle determined by a circle of the post-knot of $K$.}
\label{3handle2}
\end{figure}

\begin{figure}
\centerline{\relabelbox 
\epsfysize 6.5cm
\epsfbox{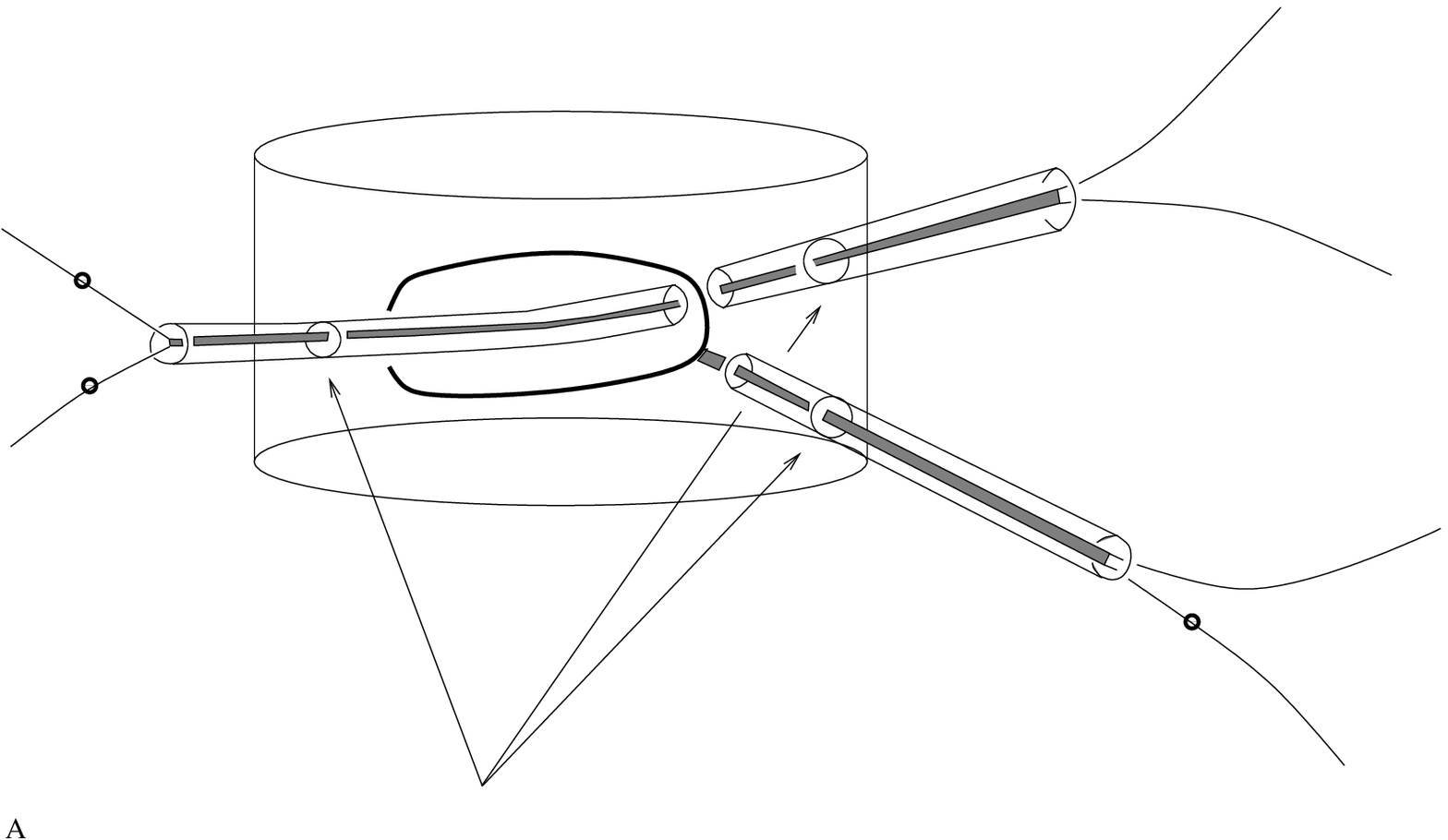}
\relabel {A}{$\s{\textrm{Intersections of the attaching sphere } S^2\textrm{  with the framed knot determined by the bands}}$}
\endrelabelbox }
\caption{Attaching sphere for a  3-handle in the Kirby diagram determined  by a knot with bands representing a knotted surface. Note that the bands appearing in the figure do not belong to the Kirby diagram, whose set of dotted circles will be the pre-knot of the knot with bands showing.}
\label{3handle}
\end{figure}

Finally we attach a 4-handle at $t=2$. By the Cellular Approximation Theorem,  this 4-handle will not affect the calculation of $\Pi_2(M,M^{(1)})=\Pi_2(M^{(3)},M^{(1)})$.

  We have thus defined a handle decomposition of the complement ${M}$ of a knotted  surface $\S$ if we are given {a knot with bands $K$ representing it}, chosen so that the the post-knot of {$K$} is a standard {diagram of the unlink}.

{Summarizing,  a Kirby diagram for {$M$}   will have a dotted circle for each circle of the pre-knot $K_-$ of $K$,} considered to be 0-framed. Then, each band of $K$ will induce a 2-handle of the complement, and the attaching region of it is determined by a framed circle $S^1 \times I$, encircling the band, with framing parallel to the core of it (thus yielding a 0-framed circle), {and going along almost the entire length of the band}{; see figure \ref{bandshandle}}. {Then we attach a 3-handle for each circle of the post-knot of $K$, which we suppose to be a standard diagram of the unlink. The attaching sphere $S^2$ of each 3-handle contains one, and only one, of the components of $K_+$ {in the region inside it}, and it intersects the framed {circles} determined by the bands of $K$ transversally, so that each connected component of the intersection (a circle $S^1$) goes around the corresponding band of $K$.}

Even though it is in not {strictly} necessary to retain the bands  {of $K$} to understand the Kirby diagram, it is useful to consider them {for determining} the fundamental crossed module of the complement.

\subsection{The Calculus}\label{calculus}

Let $\S\subset S^4=D^4_- \cup ( S^3 \times I) \cup D^4_+$ be a knotted surface defined by a knot with bands $K${, thus $\S \subset S^3 \times [0,1]$}. {Here $I=[-2,2]$.}
 We suppose that the surface $\S$, thus $K$, is oriented. Choose a regular projection $p$ of $\R^3$ onto a hyperplane of it {such} that $p$ is a regular projection  of $K$, thus defining a natural base point $*$ of $S^3$, the ``eye of the observer''; see \cite{CS}. Consider the handle decomposition of the knotted surface complement $M=S^4 \setminus \n(\S)$ just described,  and  let $M^{(1)},M^{(2)}$ and $M^{(3)}$ be the $1$-,$2$-  and $3$-{handlebodies} of it. As usual if $A\subset S^4$ is a set, we denote $A_{\leq t}=A \cap \big ( {(S^3 \times [-2, t])} \cup D^4_-\big)$.
\subsubsection{Wirtinger Relations}

The fundamental group of $M_{\leq -\e}\cong M^{(1)}$, a free group, is isomorphic with the fundamental group of the complement of the {pre-knot  $K_-$ of}   $K$ in $S^3$,  the  free group on the components of $K_-$, since $K_-$ is an unlinked union of  unknotted circles. We can define a presentation of $\pi_1(S^3  \setminus K_-)$   by considering  the Wirtinger
Presentation; see for example \cite{R}. Therefore
each arc (upper crossing) of the projection $p(K_-)$ of $K_-$ {gives a generator of} $\pi_1(S^3  \setminus K_-)$, and each crossing {yields a} relation; see figures \ref{fda} and \ref{Colour}.  For such presentation, the base point will stay at the ``eye of the observer''  of the chosen projection. 
 Notice that  we need to consider orientations on the knot diagram $p(K_-)$ 
so that these elements are well defined. Therefore, it is  at this point that
we need to introduce the (probably artificial) restriction that all the knotted surfaces that  we consider are oriented, so that we can orient the associated knot with bands $K$, providing compatible orientation of $K_-$ and $K_+$. The final result will certainly not depend on the chosen orientation of $\S$. 
\begin{figure}
\begin{center}
\includegraphics{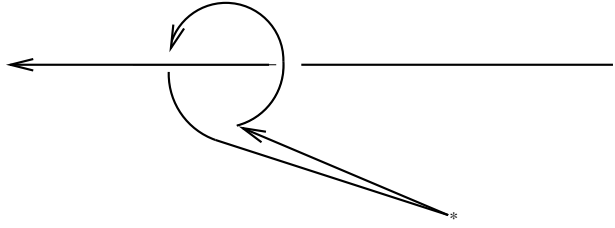}
\caption{Defining an element of the fundamental group of a knot complement. \label{fda}}
\end{center}
\end{figure}
\begin{figure}
\centerline{\relabelbox 
\epsfysize 3cm
\epsfbox{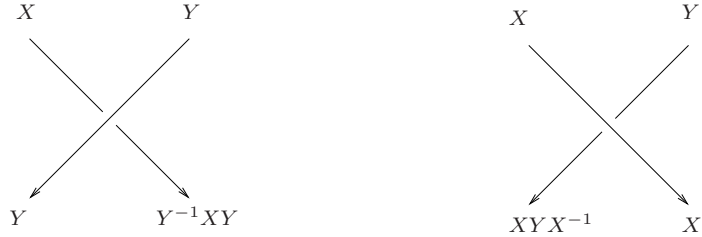} 
\relabel {X1}{$\s{X}$}
\relabel {X}{$\s{X}$}
\relabel {X3}{$\s{X}$}
\relabel {Y1}{$\s{Y}$}
\relabel {Y}{$\s{Y}$}
\relabel {Y2}{$\s{Y}$}
\relabel {Z}{$\s{Y^{-1}XY}$}
\relabel {W}{$\s{XYX^{-1}}$}
\endrelabelbox }
\caption{Wirtinger Relations.}
\label{Colour}
 \end{figure}

\subsubsection{Saddle Point Relations}

When we pass the saddle points at $t=0$, we attach $2$-handles. Therefore, by Whitehead's Theorem (Theorem \ref{Whitehead}), for  $t \in [\e,1-\e]$, the crossed module $\Pi_2(M_\l,M^{(1)}_{\l})=\Pi_2(M^{(2)},M^{(1)})$ is the free crossed module on the set of $2$-handles and their attaching maps, over the group ${\pi_1(M^{(1)})}{\cong \pi_1(S^3 \setminus K_-)}$, the free group on the set of circles of the pre-knot $K_-$  of $K$. 

Each  band of $K$ will define a 2-handle of {the knotted surface complement} $M$. However, some details are {needed} in order to specify the element of $\pi_2(M^{(2)},M^{(1)})$ defined by it as well as  its boundary in $\pi_1(M^1)$. {These are only defined up to acting and conjugation by a certain element of $\pi_1(M^{(1)})$; see  \ref{wi}.}

To make our discussion clearer, suppose that the knot with bands $K$ is such that each band of $K$  always has the same side facing upwards, for the chosen  regular projection of $K$. Standard arguments prove that any knot with bands can be isotoped so that it is in this form, whilst keeping	 the post-knot of it as being a standard diagram for the unlink. From now on any regular projection of a knot with bands will be supposed to have this special form.

To  specify the  element of $\pi_2(M^{(2)},M^{(1)})$ induced by each band, we will need to consider an orientation of the core of it (therefore defining the attaching map of the associated 2-handle up to isotopy), as well as an  arc (upper crossing) of the  band in the {chosen regular   projection $p(K)$  of $K$.} This is similar  to  the definition of the Wirtinger presentation of knot complements.

The exact {definition of these} elements of $\pi_2(M^{(2)},M^{(1)})$ is the following; cf. figure \ref{bandshandle}. {Choose a base point $*'$ on the upper part (with respect to the projection $p$) of the {framed circle} $S^1 \times I$ determined by the band.  Consider a based circle $(c,*')$ {contained in the {framed circle,}} so that $c$ goes around the band. The circle $(c,*')$ is the boundary of a certain based disk $(U,*')$ embedded in the 2-handle associated with the band, parallel to its core, by definition of attachment of 2-handles. Therefore the based disk $(U,*')$ defines an element of $\pi_2(M^{(2)}, M^ {(1)},*')$ as in \ref{wi}, and its boundary  in {$\pi_1(M^{(1)},*')$} is exactly $c$.  Note that the attaching map of the attached 2-handle is well defined since the core of the band is oriented.}

Suppose that $*'$ can be connected to the base point $*$ (the ``eye of the observer'' of $p(K)$) by a straight line which does not intersect the knot with bands. The natural {isomorphism} $\Pi_2(M^{(2)},M^{(1)},*')\to \Pi_2(M^{(2)},M^{(1)},*)$ defined by this curve determines the element of $\pi_2(M^{(2)},M^{(1)},*)$ specified by an arc  of a band in $p(K)$. (It is easy to see that this element depends only on the arc of a band  to which the base point $*'$ belongs.)

 To determine the boundary in $\pi_1(M^{(1)})$ of the elements  of $\pi_2(M^{(2)},M^{(1)})$ defined by arcs of bands,  we can use the following proposition:
\begin{Proposition}
{The boundary $\d(e)\in \pi_1(M^{(1)})\cong \pi_1(S^3 \setminus K_-)$ of the element $e\in \pi_2(M^{(2)},M^{(1)})$ determined by an arc of a band {satisfies} the relations of figure \ref{attach2}. Note that $\pi_1(M^{(1)})$ is {isomorphic} to the fundamental group of the {complement of the} pre-knot $K_-$ of $K$, itself presented by the Wirtinger Presentation of knot complements.}
\end{Proposition}
\begin{Proof}
This follows immediately  from the definition of the Wirtinger Presentation as well as the definition of the elements determined by an arc of a band. We refer to figure \ref{Boundary} for the proof in one particular  case (second case of figure \ref{attach2}). 
\end{Proof}

\begin{figure}
\centerline{\relabelbox 
\epsfysize 15cm
\epsfbox{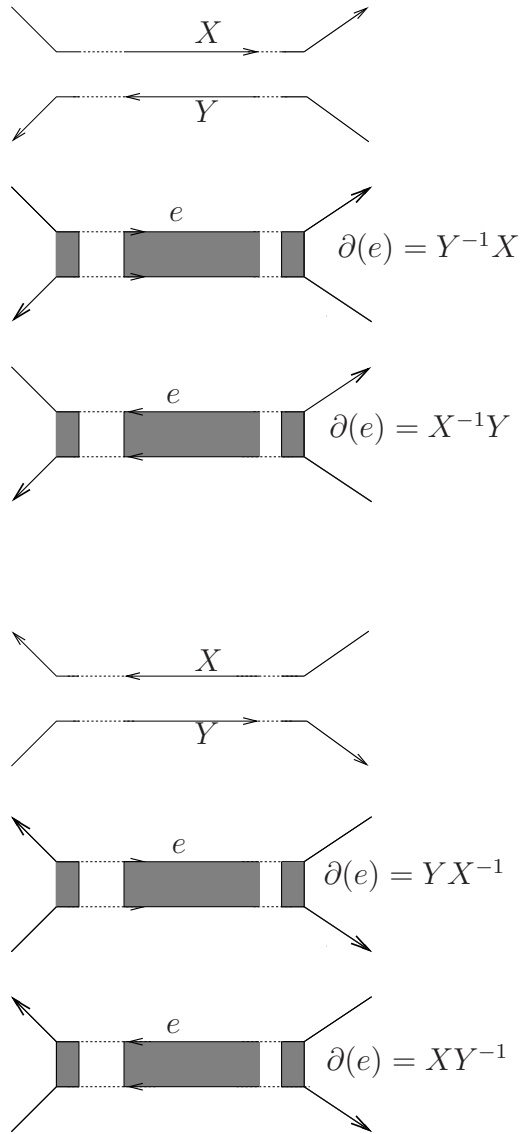} 
\relabel {X}{$X$}
\relabel {X5}{$X$}
\relabel {Y}{$Y$}
\relabel {Y5}{$Y$}
\relabel {f}{$e$}
\relabel {g}{$e$}
\relabel {e}{$e$}
\relabel {e2}{$e$}
\relabel {rel1}{$\d(e)=Y^{-1}X$}
\relabel {rel2}{$\d(e)=X^{-1}Y$}
\relabel {rel3}{$\d(e)=YX^{-1}$}
\relabel {rel4}{$\d(e)=XY^{-1}$}
\endrelabelbox }
\caption{The boundary  $\d(e) \in \pi_1(M^{(1)})$ of  the element {$e \in \pi_2(M^{(2)},M^{(1)})$ determined} by an arc of a band of a knot with bands. The knots without bands appearing in the figure  are the pre-knots of the {remaining}.}
\label{attach2}
 \end{figure}

\begin{figure}
\centerline{\relabelbox 
\epsfysize 5cm
\epsfbox{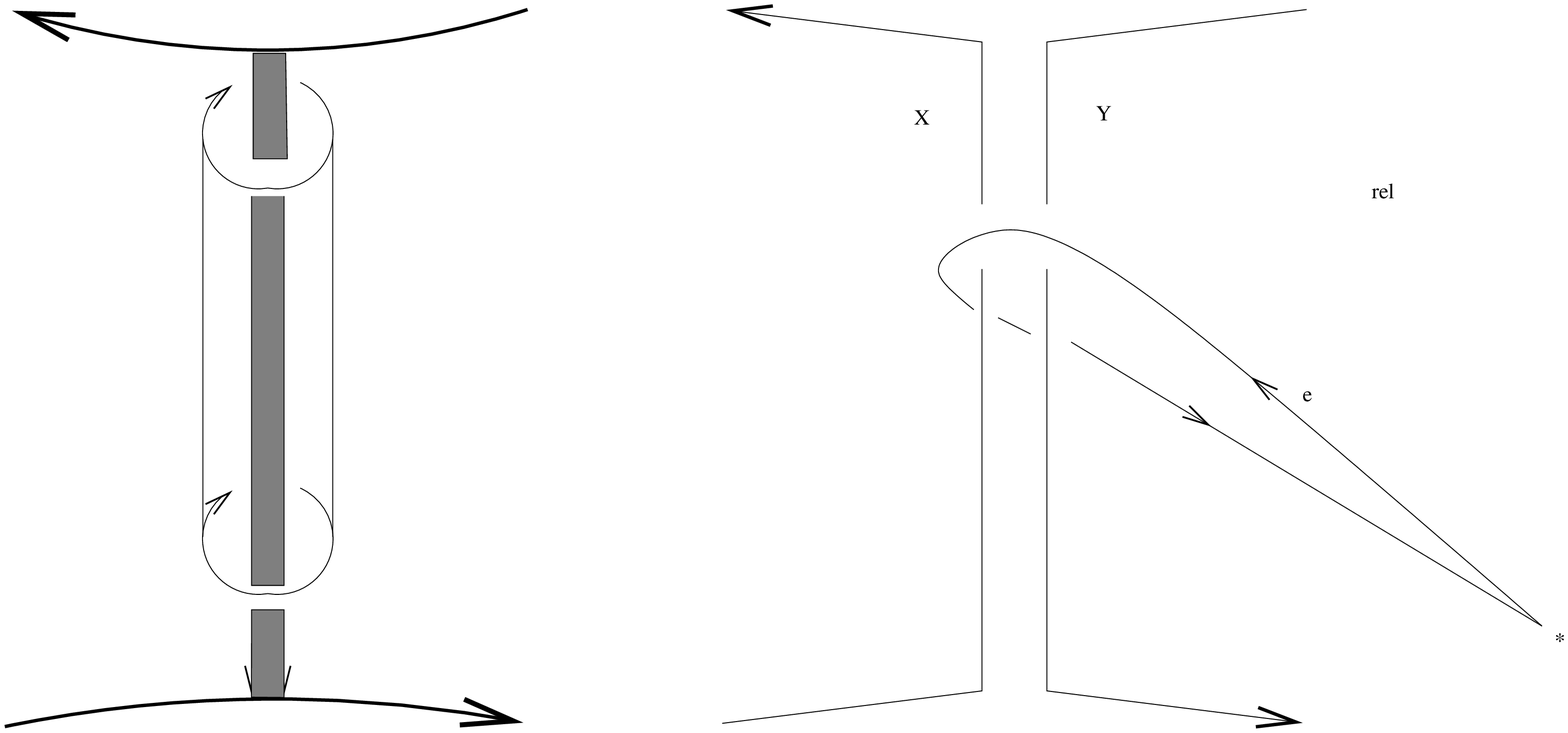} 
\relabel {X}{$Y$}
\relabel {Y}{$X$}
\relabel {e}{$\d(e)$}
\relabel {rel}{$\d(e)=Y^{-1}X$}
\relabel {*}{$*$}
\endrelabelbox }
\caption{Calculation of the boundary in $\pi_1(M^{(1)})$ of the element {$e$} of $\pi_2(M^{(2)},M^{(1)})$ determined by an arc of a band of a knot with bands.}
\label{Boundary}
 \end{figure}

The element of $\pi_2(M^{(2)},M^{(1)})$ determined by a band of a knot with bands depends on the arc chosen.  The exact dependence on the arc is described in the following proposition: 
 \begin{Proposition}
The elements of $\pi_2(M^{(2)},M^{(1)})$ determined by {arcs of  bands of a knot with bands} satisfy the relations of figure \ref{relations}.
 \begin{figure}
\centerline{\relabelbox 
\epsfysize 3cm
\epsfbox{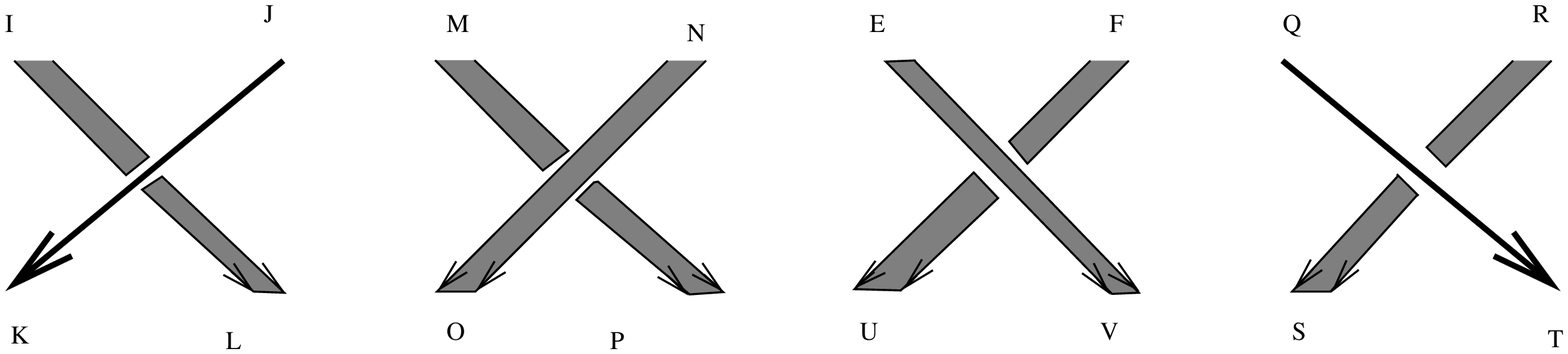} 
\relabel {E}{$\s{e}$}
\relabel {F}{$\s{f}$}
\relabel {U}{$\s{efe^{-1}}$}
\relabel {V}{$\s{e}$}
\relabel {I}{$\s{e}$}
\relabel {J}{$\s{X}$}
\relabel {K}{$\s{X}$}
\relabel {L}{$\s{X^{-1} \t e}$}
\relabel {M}{$\s{e}$}
\relabel {N}{$\s{f}$}
\relabel {O}{$\s{f}$}
\relabel {P}{$\s{f^{-1}ef}$}
\relabel {Q}{$\s{X}$}
\relabel {R}{$\s{e}$}
\relabel {S}{$\s{X\t e}$}
\relabel {T}{$\s{X}$}
\endrelabelbox }
\caption{Relations satisfied by the elements of  $\pi_2(M^{(2)},M^{(1)})$ specified by an arc of a band {of a knot with bands.}}
\label{relations}
 \end{figure}
\end{Proposition}
\begin{Proof}
This follows essentially from the fact that if we use two different paths $\g_1$ and $\g_2$ from {the base point of a CW-space $N$ to the base point of a 2-cell $c$} to identify the 2-cell {$c$} with an element of {$\pi_2(N^2,N^1,*)$,} then the corresponding   elements $c_{\g_1}$ and  $c_{\g_2}$ of {$\pi_2(N^2,N^1,*)$} are related by $c_{\g_1}=(\g_1 \g_2^{-1}) \t c_{\g_2}$.  To prove the middle bits, we also need to use the second condition of the definition of crossed modules, together with the relations of figure \ref{attach2}.
\end{Proof}

We have:

\begin{Proposition} \label{precise1} Let $\S$ be an oriented knotted surface, and let $K$ be a knot with bands representing $\S$, provided with some regular projection. As usual, suppose that the post-knot $K_+$ {of $K$} is a standard diagram for the unlink, and that each band  of $K$ always has the same side facing upwards. Let $M=S^4 \setminus \n(\S)$ be the complement of an open regular neighborhood of $\S${, with handle decomposition as above.} For each band of $K$ choose an arc of it,  therefore defining an element of $\pi_2(M^{(2)},M^{(1)})$, as well as its boundary in $\pi_1(M^{(1)})$, the free  group on the set of circles of the pre-knot  $K_-$ of $K$ (itself {presented}  by the Wirtinger Presentation of knot complements); see figure \ref{attach2}. Then  the crossed module $\Pi_2(M^{(2)},M^{(1)})$ is the free crossed module on this map {from the set of bands of $K$ into  $\pi_1(M^{(1)}){\cong \pi_1(S^3 \setminus K_-)}$.}
\end{Proposition}

\subsubsection{{The 2-Relations at Maximal Points}}

Given the drawing of the attaching {sphere} of a $3$-handle in {figure 
\ref{3handle}},  it is possible to determine the {2-relations} (implied by Theorem \ref{Attach3}) coming from the attachment  of a $3$-handle at a maximal point. First of all we  need to determine the element $r \in \pi_2(M^{(2)},M^{(1)})$  induced by the attaching map of each 3-handle; see Theorem \ref{Attach3}.   We will need the following intuitive lemma, whose straightforward  proof is left to the reader.

Consider the 2-sphere  $S^2$ having the north pole as a base point $*$. Remove {2-disks} $D^2_i$, where $i=1,2,{\ldots} ,n$ from $S^2$, whose middle points lie  in the equator of $S^2$,  with base points $*_i$ at their highest points; see figure \ref{lemma}, for the case $n=3$.  Consider the obvious paths $\g_i$ from $*_i$ to $*$ along a meridian. Then $S^2$ is obtained from $S^2 \setminus ({\cup_{i=1}^n }D^2_i)$ by attaching a 2-cell for each $i=1,{\ldots} ,n$, along the obvious attaching map $X_i$ shown in figure \ref{lemma}. The paths $*_i \to *$ permit us  to associate elements $e_i$ of $\pi_2\big (S^2,S^2 \setminus ({\cup_{i=1}^n }D^2_i)\big )$ to each disk $D^2_i$.
\begin{Lemma} \label{A}
{The product $e_1e_2{\ldots} e_n \in \pi_2\big (S^2,S^2 \setminus ({\cup_{i=1}^n }D^2_i)\big )$ coincides
with the element of $\pi_2\big (S^2,S^2 \setminus ({\cup_{i=1}^n }D^2_i)\big)$  naturally defined by the oriented $S^2$.}
\end{Lemma}

\begin{figure} 
\centerline{\relabelbox 
\epsfysize 8cm
\epsfbox{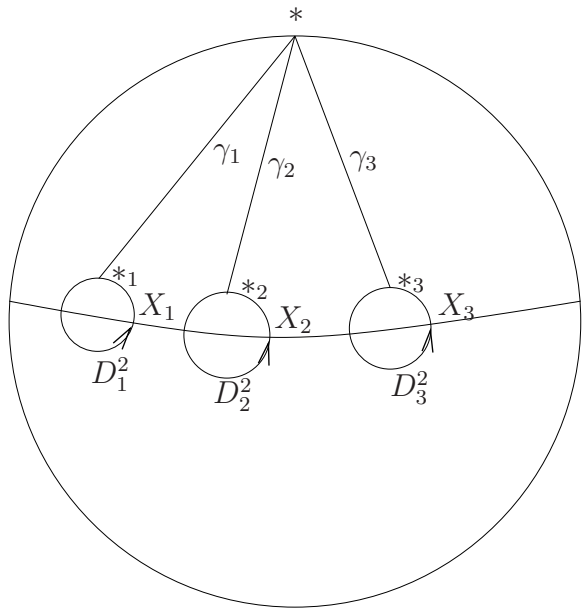}
\relabel{A}{$D^2_1$}
\relabel{B}{$D^2_2$}
\relabel{C}{$D^2_3$}
\relabel{e}{$*_1$}
\relabel{f}{$*_2$}
\relabel{g}{$*_3$}
\relabel{h}{$*$}
\relabel{x}{$\gamma_1$}
\relabel{y}{$\gamma_2$}
\relabel{z}{$\gamma_3$}
\relabel{X}{$X_1$}
\relabel{Y}{$X_2$}
\relabel{Z}{$X_3$}
\endrelabelbox }
\caption{\label{lemma} Figure relative to Lemma \ref{A}.} 
\end{figure}

Let $K$ be an oriented knot with bands so that the post- and pre-knots {$K_+$ and $K_-$} of it are oriented unlinks, with the post-knot of $K$ being a standard diagram of the unlink. Suppose also that each band of $K$ always has the same side facing upwards. Consider the natural handle decomposition {of} the complement $M$ of the oriented knotted surface $\S$ determined by $K$. Let $c_1,{\ldots} ,c_m$ be the circles of the post-knot of $K$, thus there exists one $3$-handle of $M$ for each $c_i$. For each $i=1,{\ldots} ,m$,  let $S^2_i$ be a sphere (the attaching sphere of the associated 3-handle) containing the circle $c_i$ {in the region inside it}, so that the other circles $c_j$ are outside it. {As before, suppose that $S^2_i$ intersects the framed circles determined by the bands of $K$, transversally, along the equator, so that each connected component of the intersection is a circle $S^1$ encircling the corresponding band; see Lemma \ref{B}.} Recall that this type of intersections are called essential.

  Each sphere $S^2_i$, where $i=1,2,{\ldots} ,m$, has a natural orientation induced by the orientation of a ball in $S^3$. Let ${ \{e^1_i,{\ldots} ,e^{n_i}_i\}}$ be the set of essential intersections  of the framed circles determined by the bands of  $K$ with  $S^2_i$ (note that each framed circle may intersect $S^2_i$ more than once), ordered as in figure \ref{lemma}. By {``intersections''} we mean connected components of the intersection of $S^2_i$ with the framed circles {determined by the bands}.

The core of each band is oriented, by assumption. Let $\theta^j_i=1$ if, with respect to the attaching sphere $S^2_i$, the bit of band determining  $e^j_i$ is pointing  outwards, and $-1$ otherwise.

  Each intersection $e^j_i$ induces an element {$e^j_i\in\pi_2(M^{(2)},M^{(1)})$}, determined by the arc of band which $e^j_i$ is encircling. This group element is the one provided from the fact that the intersection $e^j_i$ bounds a disk embedded in the corresponding 2-handle, parallel to its core. {This} disk is {a connected component of the} intersection of the attaching sphere  $S^2_i$ with the 2-handle; see  Lemma \ref{B}. {This discussion implies:}

\begin{Theorem}  \label{precise2}
{The crossed module $\Pi_2(M,M^{(1)})$ can be presented by the  map from the set of  bands of $K$ into  $\pi_1(S^3 \setminus K_-)$ defined in Proposition \ref{precise1}, considering a  2-relation:}

$$\prod_{j=1}^{n_i} (e^j_i)^{\theta^j_i}=1,$$
{{well defined up to cyclic permutations,}
 for each circle $c_i$ of the post-knot of $K$.}
\end{Theorem}
\begin{Proof}
The crossed module $\Pi_2(M^{(2)},M^{(1)})$ is the free crossed module on the set of bands of $K$, {and the attaching maps of the associated 2-handles}; see Proposition \ref{precise1}. The manifold $M^{(3)}$  obtained from attaching  the 3-handles {determined by the circles of $K_+$ to $M^{(2)}$} is homotopic to the space obtained by attaching a 3-cell {along} each circle $S^2_i$.

We now need to apply Theorem \ref{Attach3}. The set of {2-relations} yielding the  fundamental crossed module $\Pi_2(M^{(3)},M^{(1)})$ is given by the elements $r_i$ of $\pi_2(M^{(2)},M^{(1)})$ determined by the attaching sphere $S^2_i$  of each 3-handle; here $i=1,{\ldots} ,m$, where $m$ is the number of circles of the post-knot $K_+$ of $K$. In terms of the generators of $\pi_2(M^{(2)},M^{(1)})$ determined by the bands of the knot with bands, each of these elements is  exactly given by the formula  $r_i=\prod_{j=1}^{n_i} (e^j_i)^{\theta^j_i}$. This follows from  lemmas \ref{A} and \ref{B}.  
\end{Proof}

Note that the {2-relation} motivated by the attachment of a $3$-handle {at a maximal point} in principle depends on the chosen  sphere $S^2$ containing the corresponding circle of the post-knot $K_+$ {of $K$} {in the region inside it}. The configuration of this {attaching} sphere will usually be like in figure \ref{3handle}; so some bands of $K$ may be entirely {in the region inside}  the chosen attaching sphere{. In particular they do not appear in the {2-relation} motivated by the attachment of {the 3-handle.}} This is because the intersection of the attaching sphere with the framed {circle} determined by the band is empty; see figure \ref{3handlerel} for one such example.

\begin{figure} 
\centerline{\relabelbox 
\epsfysize 4cm
\epsfbox{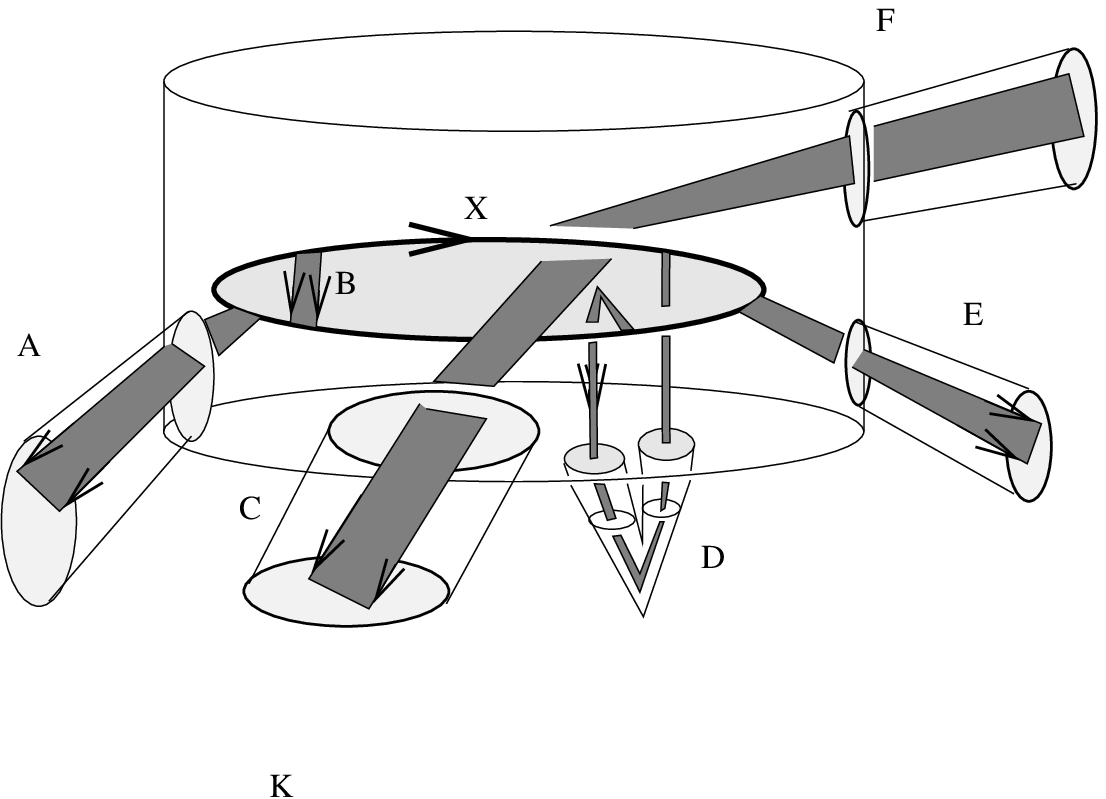}
\relabel {A}{$a$}
\relabel {B}{$b$}
\relabel {C}{$c$}
\relabel {D}{$d$}
\relabel {E}{$e$}
\relabel {F}{$X^{-1} \t c$}
\relabel {X}{$X$}
\relabel {K}{$\implies  acd d^{-1} e (X^{-1} \t c^{-1})=1$}
\endrelabelbox }
\caption{\label{3handlerel} {The 2-relations implied by the attachment of a 3-handle at a maximal point in one particular example.}} 
\end{figure}

\subsubsection{{Simple Examples}}\label{SimpleExamples}
Consider the surface $\S_1$ represented by the knot with bands ${K_1}$ of  figure \ref{ExC}. Therefore $\S_1$ is a trivial embedding of a sphere $S^2$ in $S^4$.  Let $M_1$ be the complement of an open regular neighborhood of $\S_1$, provided with the handle decomposition determined by ${K_1}$. The calculation of $ \Pi_2(M_1,M_1^{(1)})$ appears in figure \ref{ExC}. This permits us to conclude that:
$${\Pi_2(M_1,M_1^{(1)}) ={\cal U}\left (\{f\} \ra{ f \mapsto 1} F(X); X^{-1} \t f=1 \right),}$$
where $F(X)$ is the free group on the variable $X$.

\begin{figure} 
\centerline{\relabelbox 
\epsfysize 8cm
\epsfbox{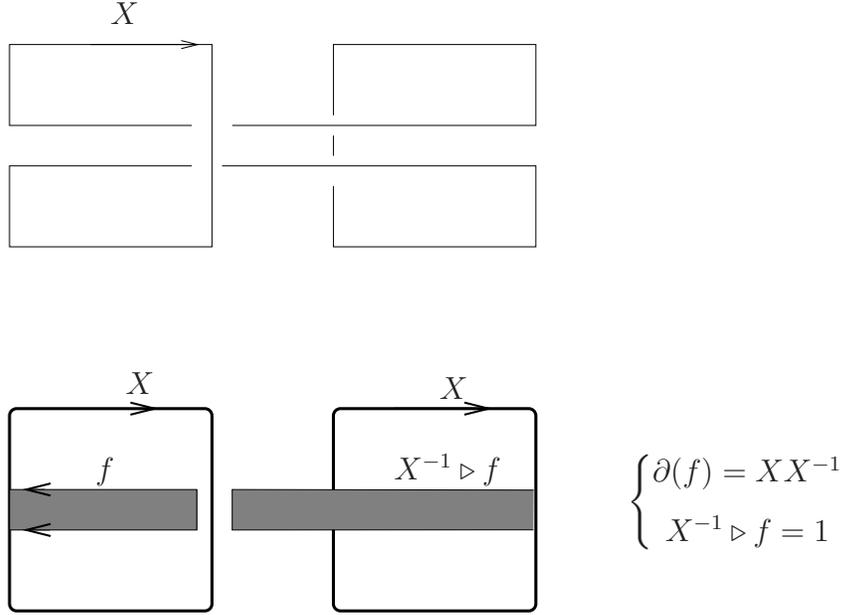}
\relabel {X}{$X$}
\relabel {Y}{$X$}
\relabel {Z}{$X$}
\relabel {e}{$X^{-1} \t f$}
\relabel {f}{$ f$}
\relabel {rel}{${{\left \{ \begin{CD} \d( f)=XX^{-1}\\ X^{-1} \t f= 1\end{CD} \right .}} $}
\endrelabelbox }
\caption{\label{ExC} {Calculation} of $\Pi_2(M_1,M^{(1)})$.  {The knot on top is the pre-knot of the knot with bands $K_1$ on the bottom.}} 
\end{figure}

Consider the knot with bands shown in figure \ref{ExA}, let $M_2$ be the complement of the knotted surface $\S_2$ represented by it. From figure \ref{ExA} we have:
$${\Pi_2(M_2,M_2^{(1)}) ={\cal U}\left (\{e\} \ra{ e \mapsto 1} F(X); e=1 \right).}$$
\begin{figure} 
\centerline{\relabelbox 
\epsfysize 6cm
\epsfbox{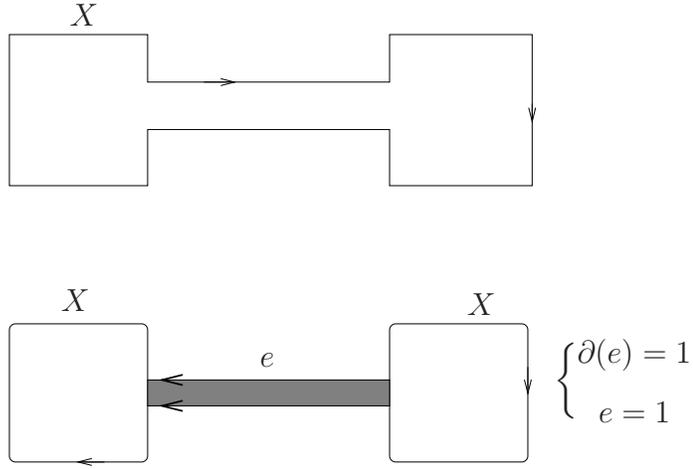}
\relabel {X}{$X$}
\relabel {X1}{$X$}
\relabel {X2}{$X$}
\relabel {e}{$e$}
\relabel {rel}{$ \left \{ \begin{CD} \d(e)=1 \\ e=1  \end{CD} \right .$}
\endrelabelbox }
\caption{ Calculation of $\Pi_2(M_2,M_2^{(1)})$. {The knot on top is the pre-knot of the knot with bands $K_2$ on the bottom.}}
\label{ExA}
\end{figure}

Likewise, consider the complement $M_3$ of the knotted surface $\S_3$ represented by the knot with bands $K_3$ of figure \ref{ExB}. Then: 
$${\Pi_2(M_3,M_3^{(1)})={\cal F}\left (\{e\} \ra{ e \mapsto X^{-1}Y} F(X,Y)\right ).}$$
Here $F(X,Y)$ is the free group on the variables $X$ and $Y$.
Note that in {this case}  we can find a sphere containing the post-knot of {$K_3$} {in the region inside it} that does not intersect the unique band of {$K_3$}, hence there are no {2-relations} motivated by the attachment of $3$-handles; see Theorem \ref{precise2}.

\begin{figure} 
\centerline{\relabelbox 
\epsfysize 6cm
\epsfbox{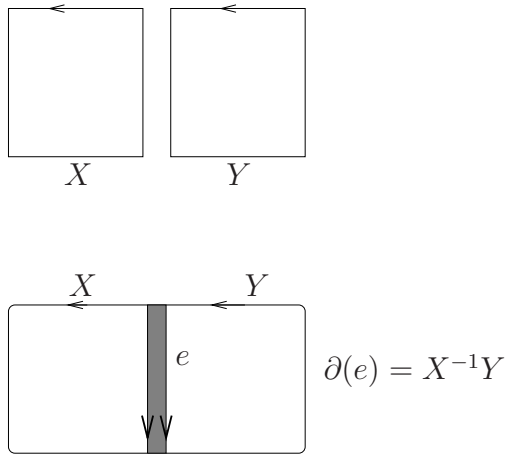}
\relabel {X}{$X$}
\relabel {Y}{$Y$}
\relabel {X2}{$X$}
\relabel {Y2}{$Y$}
\relabel {e}{$e$}
\relabel {rel}{$ \d(e)=X^{-1}Y$}
\endrelabelbox }
\caption{ Calculation of $\Pi_2(M_3,M_3^{(1)})$. {The knot on top is the pre-knot of the knot with bands $K_3$ on the bottom.}} 
\label{ExB}
\end{figure}

\subsection{Spun Trefoil}

A knot with bands {representing} the Spun Trefoil {(a knotted sphere $S^2$ embedded in $S^4$)} appears in figure \ref{STKWB}. This can be obtained for example from the marked vertex diagram of it in \cite[Figure 5]{L}, by switching the marked vertices to bands, and isotoping the final result so that the post-knot of it is a standard {diagram of the unlink.}

\begin{figure}
\centerline{\relabelbox 
\epsfysize 5cm
\epsfbox{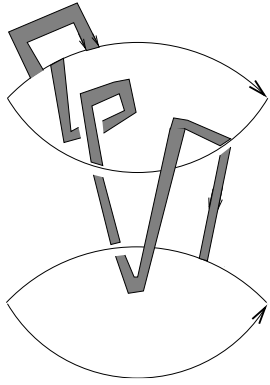}
\endrelabelbox}
\caption{Representing the Spun Trefoil by using a knot with bands.}
\label{STKWB}
\end{figure}

Let $M=S^4 \setminus \n(\S)$ be the complement of the Spun Trefoil $\S$.
We display in figure \ref{Trefoilkwb} the calculation of $\Pi_2(M,M^{(1)})$.
 \begin{figure}
\centerline{\relabelbox 
\epsfysize 17cm
\epsfbox{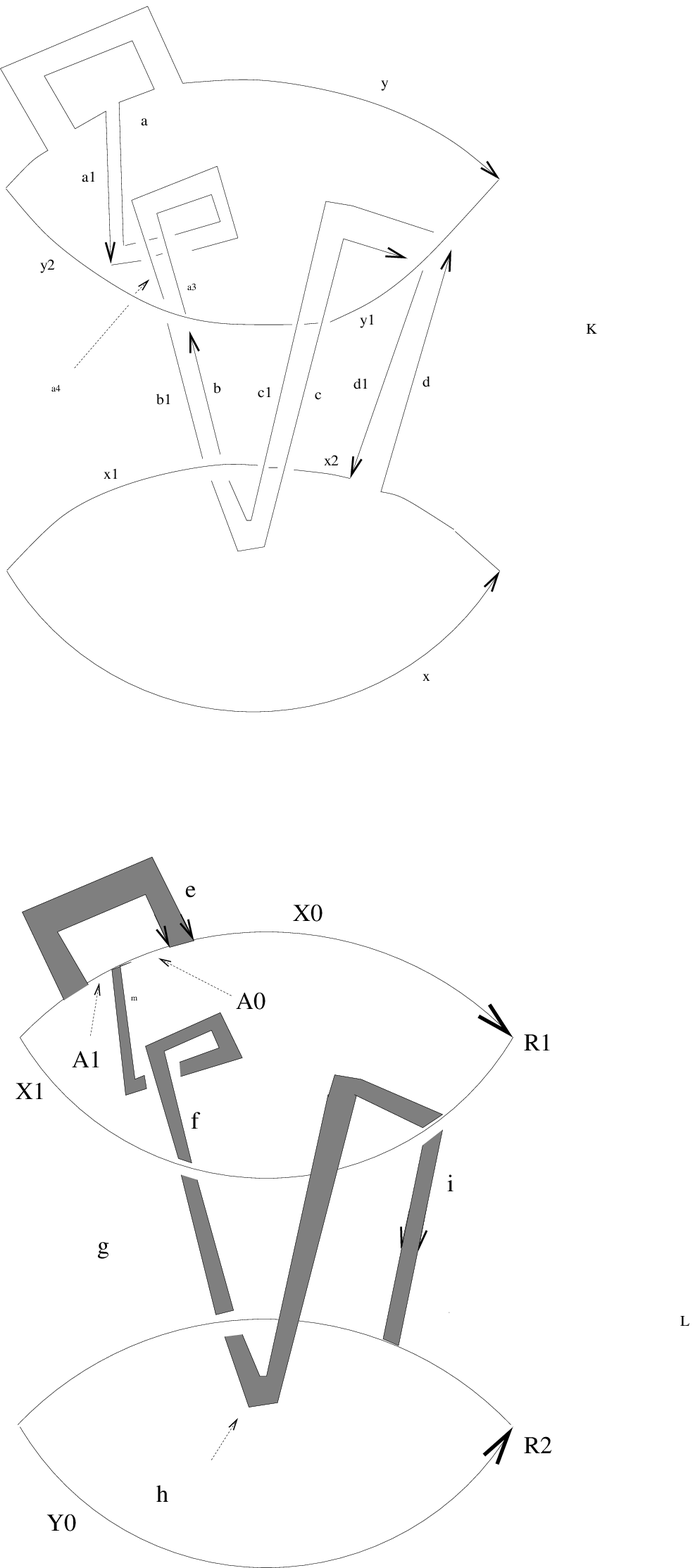}
\relabel{X0}{$\s{X}$}
\relabel{X1}{$\s{X}$}
\relabel{A0}{$\s{A}$}
\relabel{A1}{$\s{A}$}
\relabel{e}{$\s{e}$}
\relabel{f}{$\s{f}$}
\relabel{m}{$\s{f}$}
\relabel{g}{$\s{X^{-1} \t f}$}
\relabel{h}{$\s{Y^{-1}X^{-1} \t f}$}
\relabel{i}{$\s{X^{-1} Y^{-1}X^{-1} \t f}$}
\relabel{Y0}{$\s{Y}$}
\relabel{R1}{${\s{(X^{-1}\t f)( Y^{-1} X^{-1} \t f^{-1})( X^{-1}Y^{-1} X^{-1} \t f)=1}}$}
\relabel{R2}{${{\s{(X^{-1}\t f^{-1})( X^{-1}Y^{-1} X^{-1} \t f^{-1})( Y^{-1} X^{-1} \t f)=1}}}$}
\relabel{x}{$\s{Y}$}
\relabel{x1}{$\s{Y}$}
\relabel{x2}{$\s{Y}$}
\relabel{y1}{$\s{X}$}
\relabel{y2}{$\s{X}$}
\relabel{a}{$\s{A}$}
\relabel{a1}{$\s{A}$}
\relabel{a3}{$\s{A}$}
\relabel{a4}{$\s{A}$}
\relabel{b}{$\s{B}$}
\relabel{b1}{$\s{B}$}
\relabel{c}{$\s{C}$}
\relabel{c1}{$\s{C}$}
\relabel{d}{$\s{Y}$}
\relabel{d1}{$\s{Y}$}
\relabel{y}{$\s{X}$}
\relabel{K}{${\left \{ \begin{CD} C=XYX^{-1} \\ B=YXYX^{-1}Y^{-1} \\A=XYXYX^{-1}Y^{-1}X^{-1} \end{CD} \right .}$}
\relabel{L}{$\left \{ \begin{CD}  \d(f)=AA^{-1}=1 \\ \d(e)=A^{-1}X\end{CD}  \right .$}
\endrelabelbox}
\caption{Calculation of $\Pi_2(M,M^{(1)})$, where $M$ is the complement of the Spun Trefoil. }
\label{Trefoilkwb}
\end{figure}
{This permits us to conclude that $\Pi_2(M,M^{(1)})$ can be presented by the  map $\{e,f\}\to F(X,Y)$ such that $e\mapsto  A^{-1}X$, where $A=XYXYX^{-1}Y^{-1}X^{-1}$, and $f\mapsto 1$, considering also the 2-relations} $${(X^{-1}\t f)( Y^{-1} X^{-1} \t f^{-1})( X^{-1}Y^{-1} X^{-1} \t f)=1}$$ and $${(X^{-1}\t f^{-1})( X^{-1}Y^{-1} X^{-1} \t f^{-1})( Y^{-1} X^{-1} \t f)=1.}$$   Here $F(X,Y)$ is the free group on $X$ and $Y$. {Given that  $\d(f)=1$, it follows that all the elements of the form $Z \t f$, where $Z \in \pi_1(M^{(1)})$, are central in $\pi_2(M^{(2)},M^{(1)})$, by the second condition of the definition of crossed modules. In particular the two {2-relations} are equivalent, thus we can skip one of them.}

{See \ref{shl} for calculations in the Spun Hopf Link case.}

\subsubsection{The Second Homotopy  Group of the Spun Trefoil Complement}

The discussion of \ref{adec} will be needed now.  
We want to determine the kernel of the map $\pi_2(M,M^{(1)}) \to {\pi_1(M^{(1)})}$, where $M$ is the complement of the Spun Trefoil; cf. equation {(\ref{kernel})}.

 Let $K\subset S^3=\R^{3} \cup \{\infty\}$ be a knot. Suppose that  the projection on the last variable is a Morse function in $K$. Then, similarly with the 4-dimensional case, we have a handle decomposition of the complement of $K$ where  minimal points induce 1-handles of the complement and  maximal points  induce 2-handles. We have at the end to attach {one} extra 3-handle, which will cancel out  one of the 2-handles previously  attached. Therefore, the complement of the Trefoil Knot, shown in figure \ref{trefoil}, admits a handle decomposition with one 0-handle, two 1-handles and one 2-handle; see \cite{GS}, exercise $6.2.2$. Note that 2-dimensional CW-complexes with a unique $2$-cell are classified by their fundamental group, up to homotopy equivalence; see \cite{J}.
\begin{figure}
\centerline{\relabelbox 
\epsfysize 3 cm
\epsfbox{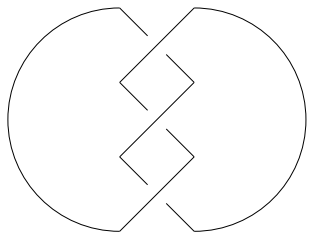}
\endrelabelbox}
\caption{Trefoil Knot.}
\label{trefoil}
\end{figure}

Consider the CW-complex $N$ with one 0-cell $\{*\}$, two 1-cells $X$ and $Y$ and a 2-cell $e$ attaching along $A^{-1}X$. Here $A=XYXYX^{-1}Y^{-1}X^{-1}$. Then $N$ is homotopic to the Trefoil Knot complement. We can prove this for example from the fact {that $N$} has a unique 2-cell and its fundamental group is isomorphic to the fundamental {group} of the Trefoil Knot  complement. In particular $\pi_2(N,*)=\{0\}$, by the well known theorem (due to  Papakyriakopoulos)  asserting that 3-dimensional (one component) knot  complements are aspherical; see \cite{Pa,R}. {On the other hand,}  we can represent $\pi_2(N,*)$ as {$ \ker \{\d\colon  \pi_2(N,N^1,*) \to \pi_1(N^1,*)\}$}. Note that {$\Pi_2(N,N^1,*)$} is the free crossed module on the map $e \mapsto A^{-1}X\in F(X,Y)$. 

Let  $M$ be the complement of the Spun Trefoil{, with handle decomposition as above}. Let ${R'}$ be the subgroup of $\pi_2(M^{(2)},M^{(1)})$ generated by the elements $B \t f$, where $B \in F(X,Y)$. Since $\d(f)=1$,  then (see \ref{adec}) the group $\pi_2(M^{(2)},M^{(1)})$ is the direct sum of {$\pi_2(N,N^{1},*)$} and ${R'}$. Moreover, ${R'}$ is the free abelian module with basis $\{f\}$ over {$\pi_1(M)={{\rm coker} (\d)\doteq G}$}, where {$\d\colon \pi_2(N,N^1,*) \to \pi_1(N^1,*)=F(X,Y)$}, with {the} obvious action of $F(X,Y)$. {In addition,}  $G=F(X,Y)/\left <XA^{-1}=1\right >$, which is isomorphic with the fundamental group of the complement of the Trefoil Knot.

Now,  {$-(X^{-1}\t f)+( Y^{-1} X^{-1} \t f)- ( X^{-1}Y^{-1} X^{-1} \t f) \in {R'}$ (note that we switched to additive notation).} Consider the abelian $G$-module:
$${R=\frac{R'}{\left < -(X^{-1}\t f)+( Y^{-1} X^{-1} \t f)-( X^{-1}Y^{-1} X^{-1} \t f)=0\right>}.}$$
Then it follows that $\pi_2(M,M^{(1)})$, as an  $F(X,Y)$-module, is the direct sum of $R$ and {$\pi_2(N,N^1,*)$.}

The CW-complex $N$ is aspherical, thus the kernel of the boundary map of {$\Pi_2(M,M^{(1)}) $} restricted to {$\pi_2(N,N^1,*)$} is trivial. 
 In particular, $\pi_2(M)=\ker \{\d\colon  \pi_2(M,M^{(1)}) \to  \pi_1(M^{(1)})\}={R}$. It follows that {the second homotopy group of the complement of the Spun Trefoil} is the free abelian module over $G$ with one generator $f$ and the relation {
$-(X^{-1}\t f)+( Y^{-1} X^{-1} \t f)-( X^{-1}Y^{-1} X^{-1} \t f)=0.$}
This is {consistent} with the calculation in \cite{L}.

\section{The ``Crossed Module Invariant'' of Homotopy Types}

\subsection{Definition of the Invariant}\label{CMI}
Let $M$ be a CW-complex. As we have seen in \ref{mainset}, the fundamental crossed module $\Pi_2(M,M^1)$ does not depend on the homotopy type of $M$, as a space, up to free products (in the category of crossed modules) with $\Pi_2(D^2,S^1)$. Given a cell decomposition of $M$, it is in principle possible to obtain a {presentation of $\Pi_2(M,M^1)$.} However, it may be difficult to distinguish  between two crossed modules (up to free products with $\Pi_2(D^2,S^1)$) presented this way.

It is possible in some cases to  distinguish between finitely generated groups  by using the Alexander's Invariant (see \cite{CF}), or by counting the number of morphisms from them into a finite group. The latter  method can also be used, with due adaptations,  in the case of finitely generated crossed modules. 
We have:

\begin{Theorem}\label{main0}
Let $M$ be a finite CW-complex with a unique $0$-cell, which we take to be its base point $*$. Let $\G=(G,E,\d,\t)$ be a finite crossed module. {Here $\d\colon E \to G$ is a group morphism and $\t$ is a left action of $G$ on $E$ by automorphisms.} The quantity:

$$I_\G(M)=\frac{\# \Hom\big (\Pi_2(M,M^1), \G\big)}{(\# E)^{b_1(M^1)}}$$
is finite, does not depend on the CW-decomposition of $M$ and it is a homotopy invariant of $M$, as a space. {Here $b_1(M^1)$ denotes the first Betti number of the 1-skeleton $M^1$ of $M$.}
\end{Theorem}
We call this homotopy invariant of  CW-complexes the ``Crossed Module Invariant''.

\begin{Proof}
The finiteness of $I_\G(M)$ follows from the fact (Theorem \ref{Attach3}), that if $M$ is finite then  the group $\pi_2(M,M^1)$ is finitely generated as a module over $\pi_1(M^1)$, itself a finitely generated group. 

 {Let} $N$ be a CW-complex with a unique $0$-cell $*'$, such that $N$ is  homotopic to $M$ as a space. By the discussion in \ref{mainset}, there exist positive integers $m$ and $n$ such that:
$$\Pi_2(M,M^1,*) \vee \Pi_2(D^2,S^1,*)^{\vee n} \cong \Pi_2(N,N^1,*')\vee \Pi_2(D^2,S^1,*')^{\vee m}.$$ 
By using the universal property defining free products of crossed modules (Example \ref{freeprod}) it follows that:
\begin{multline*}\# \Hom(\Pi_2(M,M^1,*),\G)\times \left (\# \Hom (\Pi_2(D^2,S^1,*),\G) \right) ^{n}\\ = \#   \Hom(\Pi_2(N,N^1,*),\G)\times  \left (\# \Hom (\Pi_2(D^2,S^1,*),\G) \right)^m.  
\end{multline*}  
Since  $\Pi_2(D^2,S^1,*)=(\Z,\Z,\id,\t)$, where $\t$ is the trivial action, we have that  $\#\Hom (\Pi_2(D ^2,S^1,*),\G)=\# E$. 
The result follows from the fact that we necessarily have $b_1(M^1)+n=b_1(N^1)+m$.
\end{Proof} 

Let $M$ be a compact manifold with a handle decomposition. Then the crossed module invariant $I_\G(M)$ can be calculated by  using   the crossed module $\Pi_2(M,M^{(1)})$, where as before $M^{(1)}$ is the 1-handlebody of $M$, {made out of} the {$0$-} and $1$-handles of $M$. This follows immediately from the discussion in the beginning of the previous chapter.

The Crossed Module  Invariant  {$I_\G$} has the following geometric interpretation; see \cite{FM2,FMP}. Let $\G$ be a finite crossed module. Then $\G$ has a classifying space $B_\G$ with a natural base point $*$; see \cite{BH2,B}. Consider the space $TOP\big ((M,*),(B_\G,*)\big )$  of continuous  maps $(M,*) \to (B_ \G,*)$, provided  with the $k$-ification of the compact-open topology. Then:

$$I_\G(M)= \sum_{f \in \big [(M,*), (B_\G,*)\big ]} \frac{1}{\# \pi_1\Big(TOP\big ((M,*),(B_\G,*)\big ),f\Big ) }.$$
Here $\big[(M,*), (B_\G,*)\big]=\pi_0 \Big (TOP\big ((M,*),(B_\G,*)\big)\Big)$ denotes the set of homotopy classes of maps $(M,*) \to  (B_\G,*)$.

{The results in this subsection extend in a natural way to crossed complexes; see \cite{FM2}.}
\subsubsection{Relation with Algebraic 2-Types}\label{2Types}
Recall that a 2-type is a path-connected topological space $M$ such that $\pi_k(M)=\{0\},\forall k>2$. If $M$ is a connected CW-complex with a base point $*$ which is a 0-cell, then ${T_2}(M)$ is the based cellular space defined from $M^3$ by killing all the homotopy groups $\pi_k(M)$ of $M$ with $k>2$ in the usual way; see for example \cite[{Example 4.17}]{AH}. It is well known that ${T_2}(M)$ does not depend on the CW-decomposition of  $M$ up to homotopy equivalence.   
If $M$ is a  connected CW-complex, then the CW-complex ${T_2}(M)$ is called the topological 2-type of $M$, or, more commonly, the second Postnikov section of $M$.

Let $M$ be a  path-connected topological space. Then the algebraic 2-type of $M$ is given by the triple 
$${A_2}(M)=\big (\pi_1(M),\pi_2(M),k(M)\big),$$
 where $k\in H^3(\pi_1(M),\pi_2(M))$ is the  $k$-invariant (or first Postnikov invariant) of $M$; see \cite{ML,EML,MLW}. If $M$ and $N$  are CW-complexes then their topological 2-types (up to homotopy) coincide if and only if their algebraic 2-types coincide (up to isomorphism); see \cite{Bau,EML}. 
 On the other hand we have that  $\Pi_2(M,M^1)=\Pi_2(M^3,M^1)=\Pi_2\big ({T_2}(M),{T_2}(M)^1\big )$, {where both  equalities follow from the Cellular Approximation Theorem.} Therefore, {we have:}
\begin{Theorem}
Let $\G$ be a finite crossed module. The homotopy invariant $I_\G(M)$,  depends only  on the (algebraic or topological) 2-type  of $M$. 
\end{Theorem}

Therefore a natural issue is how useful  the Crossed Module Invariant  is for separating 2-types.

\subsection{Applications to Knotted Surfaces}\label{DDD}

Theorem \ref{main0} tells us that we can {define} invariants of  knotted surfaces by considering the invariants $I_\G$ of {their complements} in $S^4$, where $\G$ is a finite crossed module. The previous chapter provides us  with an algorithm for {this type of  calculations.} 

{
Recall the construction of a  crossed module presented by a map,  with 2-relations, in \ref{general}. It follows immediately  that:}
\begin{Lemma} \label{REFER}
{Suppose that  $G$ is a free group, say on the set $L$, thus if $G'$ is a group then any map $\f_0{\colon}L\to G'$ extends uniquely to a group morphism $\f\colon G \to G'$. Let $K=\{m_1,{\ldots} ,m_k\}$ be a set provided with a map $\d_0\colon  K \to G$. Let also $\{r_1,{\ldots} ,r_n\}$ be a set of 2-relations {in} the free crossed module on the map $\d_0\colon  K \to G$.  Let $\G'=(G',E',\d',\t')$ be a crossed module. There exists a one-to-one correspondence between crossed module maps:}
 $${{\cal U}\left (K \ra{\d_0} G;r_1,{\ldots} ,r_n {=1}\right)\to \G',}$$
 and pairs of maps $\f_0\colon L\to G'$ and $\p_0\colon K \to E'$ verifying: 
\begin{enumerate}
\item  $(\f \circ \d_0)(m_i)=(\d' \circ \p_0)(m_i);i=1,{\ldots} ,k$,
\item $\psi(r_i)={1_{E'}},i=1,{\ldots} ,n.$
\end{enumerate}
Here $(\f,\p)$ is the map ${\cal F}\left ( \d_0\colon K \to  G\right ) \to \G'$ determined by $\f_0$ and $\psi_0$. In particular if 
$$r=(X_1, m_1 )^{\theta_1}(X_2, m_2 )^{\theta_2}{\ldots} (X_l,m_l)^{\theta_l},$$
{where $X_i \in G,m_i \in K$ and $\theta_i \in \Z$, for $i=1,2,\ldots l$,} then $$\psi(r)=(\f(X_1)\t' \psi_0( m_1) )^{\theta_1}(\f(X_2)\t' \psi_0( m_2) )^{\theta_2}{\ldots} (\f(X_l)  \t' \psi_0( m_l))^{\theta_l}{\in E'.}$$
\end{Lemma}

As an example, consider the knotted surface complements $M_1$, $M_2$ and $M_3$ of \ref{SimpleExamples}. By the previous lemma it follows that:
$$I_\G(M_1)=\frac{\#\{X \in G, f \in E{|}\d(f)=1, X^{-1} \t f=1\}}{\# E}=\frac{\# G}{ \# E},$$
$$I_\G(M_2)=\frac{\#\{X \in G, e \in E{|}\d(e)=1,e=1\}}{\# E}=\frac{\# G}{\# E},$$
$$I_\G(M_3)=\frac{\#\{X,Y \in G, e \in E{|}\d(e)=X^{-1}Y\}}{\# E^2}=\frac{\# G}{ \# E}.$$
In fact all these knotted surfaces are isotopic to the trivial knotted {sphere}.

Suppose that $\G=(G,E,\d,\t)$ is a finite crossed module with $G$ abelian, and $\d=1_G$, from which it follows that $E$ is abelian. Then if $M$ is the complement of the Spun Trefoil it follows that: 
$$I_\G(M)=\frac{\#\big \{(X,f) \in G \times E\colon ( X\t f)( X^2\t f^{-1 })( X^{3} \t f)\big \}}{\#E}.$$
{This agrees  with the calculation in \cite{FM1}. }

This information  {is sufficient for proving} that the Spun Trefoil is knotted. For example, consider the crossed module $\A=(\Z_2,\Z_3,\d,\t)$, where $\Z_2=\{+1,-1\}$, with the trivial boundary map $\Z_3 \to \Z_2$. The action of $\Z_2=\{+1,-1;\times\}$ in $\Z_3=\{[0],[1],[2];+\}$ is $a\t v=av$. Then $I_\A(M)=4/3$. However, if $M_1$ is the complement of the {trivial} knotted {sphere} then $I_\A(M_1)=2/3$. 

{Notice that we would have not been able to prove this (known) fact if we had used
 only the fundamental group of the complement of the  Spun Trefoil, and counting the number of
 morphisms from  it into an abelian group, since the first homology group of {a knotted surface complement} depends only on the intrinsic topology of the knotted
 surface, and not on the embedding.} {See also Remark \ref{BBB} and \ref{AAAA}.}

\subsubsection{A General Algorithm}

  Let $\G=(G,E,\d,\t)$  be a finite crossed module and let $M=S^4 \setminus \n(\S)$ be the complement of the knotted surface $\S$. To calculate $I_\G(M)$ we do not  need to determine $\Pi_2(M,M^{(1)})$ fully, which, since $\G$ verifies relations of its own, can be much more complicated than calculating $I_\G(M)$ alone.   To this end we define:

\begin{Definition}{\bf ($\G$-coloring)}
Let $K$ be a knot with bands representing some oriented knotted surface $\S$. As usual, we suppose that $K$ is provided with a regular projection; as before such that the post-knot of $K$ is a standard diagram of the unlink, and such that each band of $K$ always has the same side facing {upwards}. Let $\G=(G,E,\d,\t)$ be a finite crossed module. A $\G$-coloring of $K$ is an assignment of an element of $G$ to each arc of the {pre-knot} $K_-$ of $K$ and an element of $E$ to each arc of the bands of $K$, verifying the conditions of figures \ref{Colour}, \ref{attach2}, \ref{relations} and \ref{3handlerel}. {These} last should be interpreted in light of Theorem \ref{precise2}.
\end{Definition}

Note that the additional relations obtained when a thin component of $K$  passes under a fat component are dealt with by the Wirtinger relations for $K_-$.

By  Theorem \ref{REFER} it follows that:
\begin{Theorem}\label{genrelations}
{Let $\S \subset S^4$ be an oriented knotted surface, and choose a knot with bands $K$ representing $\S$, as well as a regular projection of $K$}. Suppose that the post-knot of $K$ is a standard diagram of the unlink{, and} {that} each band of $K$ always has the same side facing {upwards}. Let $\G=(G,E,\d,\t)$ be a finite crossed module. We have
$$I_\G\big (S^{4} \setminus \n(\S)\big)=\frac{\# {\{\G\textrm{-colorings of } K\}}}{\# E^{\# \{\textrm{circles of } K_-\}}}.$$
\end{Theorem}
{Note that, given a knot with bands $K$ representing the knotted surface $\S\subset S^4$, the handle decomposition of the complement of $\S$ in $S^4$ constructed from $K$ has a unique 0-handle and a 1-handle for each circle of the pre-knot $K_-$ of $K$.}

\begin{Remark}\label{BBB}
We prove in \cite{FMK} that the {Crossed Module Invariant} $I_\G$ is strong enough to distinguish between diffeomorphic knotted surfaces with the same fundamental group of the complement, at least in some {particular} cases; {see also \ref{AAAA}}.  In \cite{PS}, it was asserted  the existence of pairs $\S,\S'\subset S^4$ of  knotted spheres such that {$\pi_1\big(S^4\setminus \n(\S)\big)\cong\pi_1\big(S^4 \setminus \n(\S')\big)$ and  $\pi_2\big(S^4\setminus \n(\S)\big)\cong\pi_2\big(S^4 \setminus \n(\S')\big)$,} as modules over $\pi_1\big ( S^4\setminus \n(\S)\big)=\pi_1\big(S^4 \setminus \n(\S')\big)$,  but with $k\big(S^4 \setminus \n(\S) \big ) \neq k\big(S^4 \setminus \n(\S')\big)$. Since the invariant $I_\G(M)$ depends only on the algebraic $2$-type of $M$ (see \ref{2Types}), it would be interesting to determine whether $I_\G$ is strong, and practical,  enough to distinguish between pairs of embedded spheres with this property, completing the results of \cite{FMK}.
\end{Remark}

\subsubsection{Spun Hopf Link}\label{shl}

 Consider the Spun Hopf Link $\S$ obtained by spinning the Hopf Link depicted in figure \ref{hopflink}. Therefore $\S$ is an embedding of a disjoint union of {two} tori {$S^1 \times S^1$} into $S^4$. A knot with bands $K$ representing the Spun Hopf Link appears in figure \ref{movie2}. 

Let us calculate $\Pi_2(M,M^{(1)})$, where $M$ is the complement of the Spun Hopf Link, {with the}  handle decomposition determined by $K$. This calculation appears in figure \ref{Calc1}. This permits us to conclude that  $\Pi_2(M,M^{(1)})$  is the {crossed module presented by the map $\{e,f,g,h\} \to F(X,Y)$, where {$f,h \mapsto 1$}, $e \mapsto X^{-1}YXY^{-1}$ and $g \mapsto X^{-1}Y^{-1}XY$, considering the {2-relations}} $$f(X \t f^{-1}) h (Y \t h^{-1}) =1$$ and 
$$ (X \t f) f^{-1}  (Y \t h) h^{-1} =1.$$ 
 As usual, $F(X,Y)$ is the free group on the variables $X$ and $Y$.
{These {2-relations}} are equivalent, {thus} we can skip one of them. {To prove this we need to use the fact that $\d(h),\d(f)=1$, which implies that both $h$ and $f$ are central in $\pi_2(M^{(2)},M^{(1)})$.}

Let $\G=(G,E,\d,\t)$ be a finite crossed module. We can easily  calculate $I_\G\big(S^4\setminus \n(\S)\big)$, {for $\S$ the Spun Hopf Link},
 by using this {calculation}. Suppose that $\G=(G,E,\d,\t)$ is a finite crossed module with $G$ abelian and  $\d=1_G$. We have:
\begin{equation*}
I_\G(S^4\setminus {\n(\S)}){=\#\big \{X,Y \in G; f,h \in E\left|  f(X \t f^{-1}) h (Y \t h^{-1}) \right .  =1_E \big \} } .
\end{equation*}
{This particularizes to {$I_\A(S^4\setminus \n(\S))=18$,} for the case when  $\A=(\Z_2,\Z_3,\d,\t)$ is the crossed module defined above.}

\begin{figure}
\begin{center}
\includegraphics{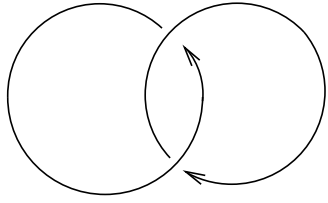}
\end{center}
\caption{Hopf Link.}
\label{hopflink}
\end{figure}
\begin{figure}
\begin{center}
\includegraphics{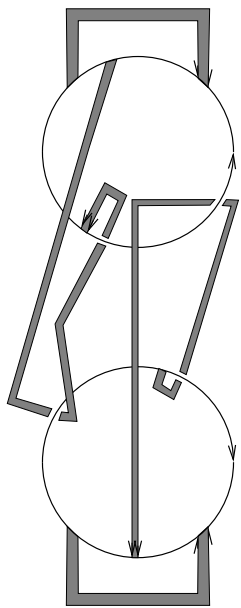}
\end{center}
\caption{A knot with bands representing the Spun Hopf Link.}
\label{movie2}
\end{figure}
\begin{figure} 
\centerline{\relabelbox 
\epsfysize 15cm
\epsfbox{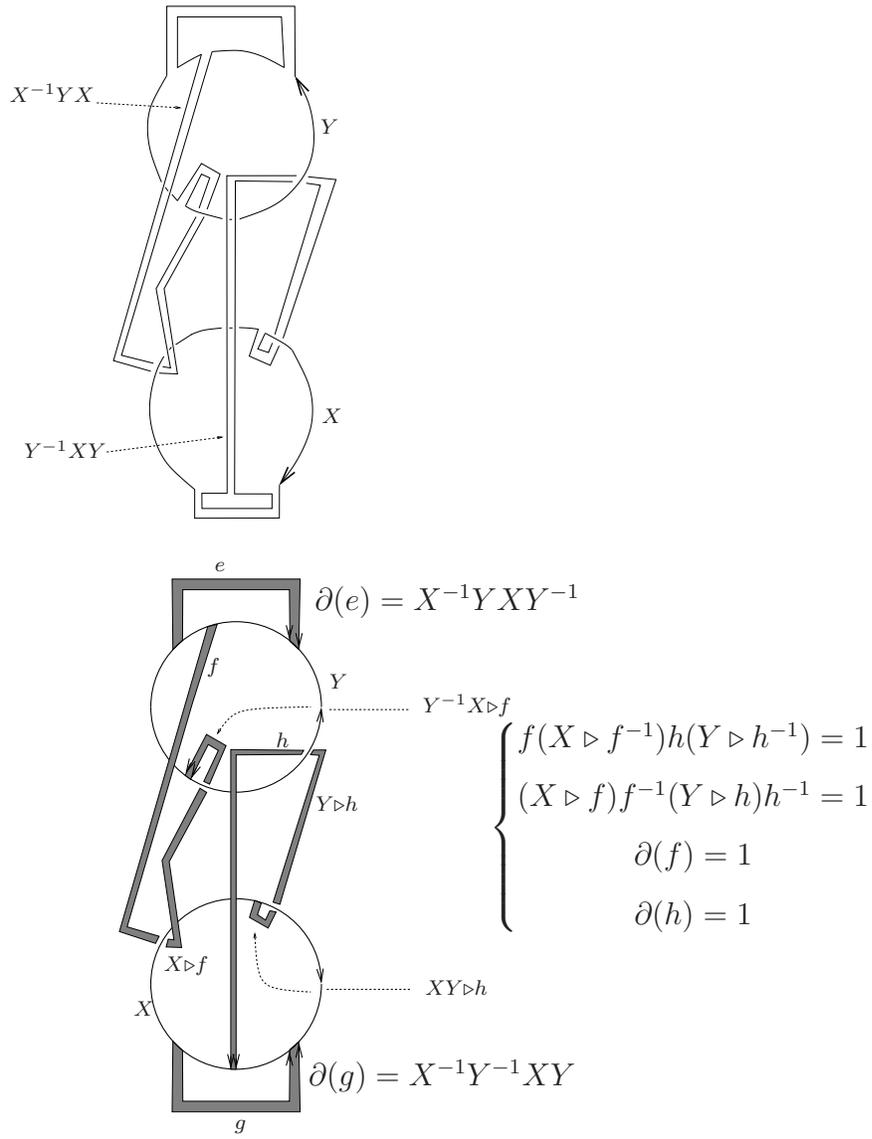}
\relabel{I}{$\s{Y}$}
\relabel{J}{$\s{X}$}
\relabel {e}{$\s{e}$}
\relabel {f}{$\s{f}$}
\relabel {g}{$\s{Y^{-1}X \t f}$}
\relabel {h}{$\s{h}$}
\relabel {i}{$\s{Y \t h}$}
\relabel {j}{$\s{g}$}
\relabel {k}{$\s{X\t f}$} 
\relabel {X}{$\s{X}$}
\relabel {Y}{$\s{Y}$}
\relabel {m}{$\s{XY \t h }$} 
\relabel {A}{$\s{X^{-1}YX}$}
\relabel {B}{$\s{Y^{-1}XY }$}
\relabel {rel}{{${\left \{ \begin{CD}  {f(X \t f^{-1}) h (Y \t h^{-1})} =1 \\ { (X \t f) f^{-1}  (Y \t h) h^{-1} }=1\\ \d(f)=1\\ \d(h)=1 \end{CD} \right . }$}}
\relabel {aa}{$ \d(e)=X^{-1}YXY^{-1}$}
\relabel {bb}{$ \d(g)=X^{-1}Y^{-1}XY$}
\endrelabelbox }
\caption{\label{Calc1} Calculation of $\Pi_2(M,M^{(1)})$ where $M=S^4 \setminus \n(\S)$. Here $\S$ is the Spun Hopf Link.}
\end{figure}

If $T$ is the trivial embedding of a disjoint union of two tori we have, for {any finite}  crossed module $\G=(G,E,\d,\t)$:

$$I_\G(S^4 \setminus \n(T))=\left (\frac{(\#G) (\# \ker \d)^2}{\# E}\right )^2,$$
which specializes to: 
$$I_\G(S^4 \setminus \n(T))=(\#G)^2(\#E)^2, $$
in the particular case for which $G$ is abelian and $\d=1_G$.
{Comparing with the value for {$I_\A(S^4 \setminus \n(\S))$} when $\S$ is the Spun Hopf Link, proves that {$\Sigma$} is  knotted.}

We can also determine the second homotopy group of the {complement $M$ of the} Spun Hopf Link from the presentation  of the fundamental crossed module $\Pi_2(M,M^{(1)})$. Proceeding as in the case of the Spun Trefoil, it follows {that $\pi_2(M)$}  is the {quotient of the free abelian module over $G=\{X,Y| XY=YX\}$ {(the fundamental group of $M$)}  with generators $\{z,f,h\}$ by the relation $X \t f -f + (Y \t h) -h =0.$}

To prove this, we need to use the fact that the CW-complex $N$ constructed with two 1-cells $X$ and $Y$ and two 2-cells $e$ and $g$, attaching in a way such that $\d(e)= X^{-1}YXY^{-1}$, whereas $\d(g)= X^{-1}Y^{-1}XY$, is such that $\pi_1(N)=G$ and $\pi_2(N)$ is the free {abelian module over $\Z[G]$} with basis $z$. This follows from the fact that $N\cong T^2 \vee S^2$, easy to prove. {Here $T^2=S^1 \times S^1$ is the torus.}

\subsubsection{Final Example}\label{AAAA}

\begin{figure}
\begin{center}
\includegraphics{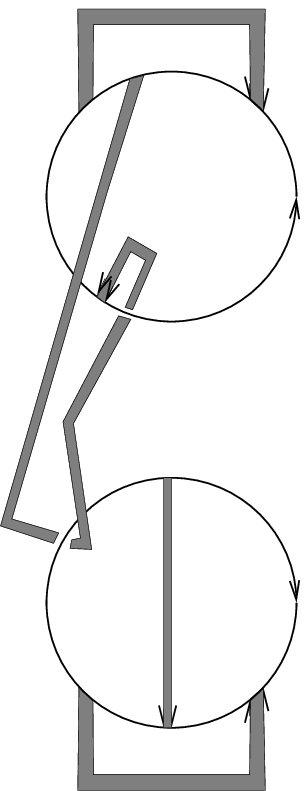}
\end{center}
\caption{A knot with bands $K'$ representing the knotted surface $\S'$.}
\label{movie2a}
\end{figure}
\begin{figure} 
\centerline{\relabelbox 
\epsfysize 7cm
\epsfbox{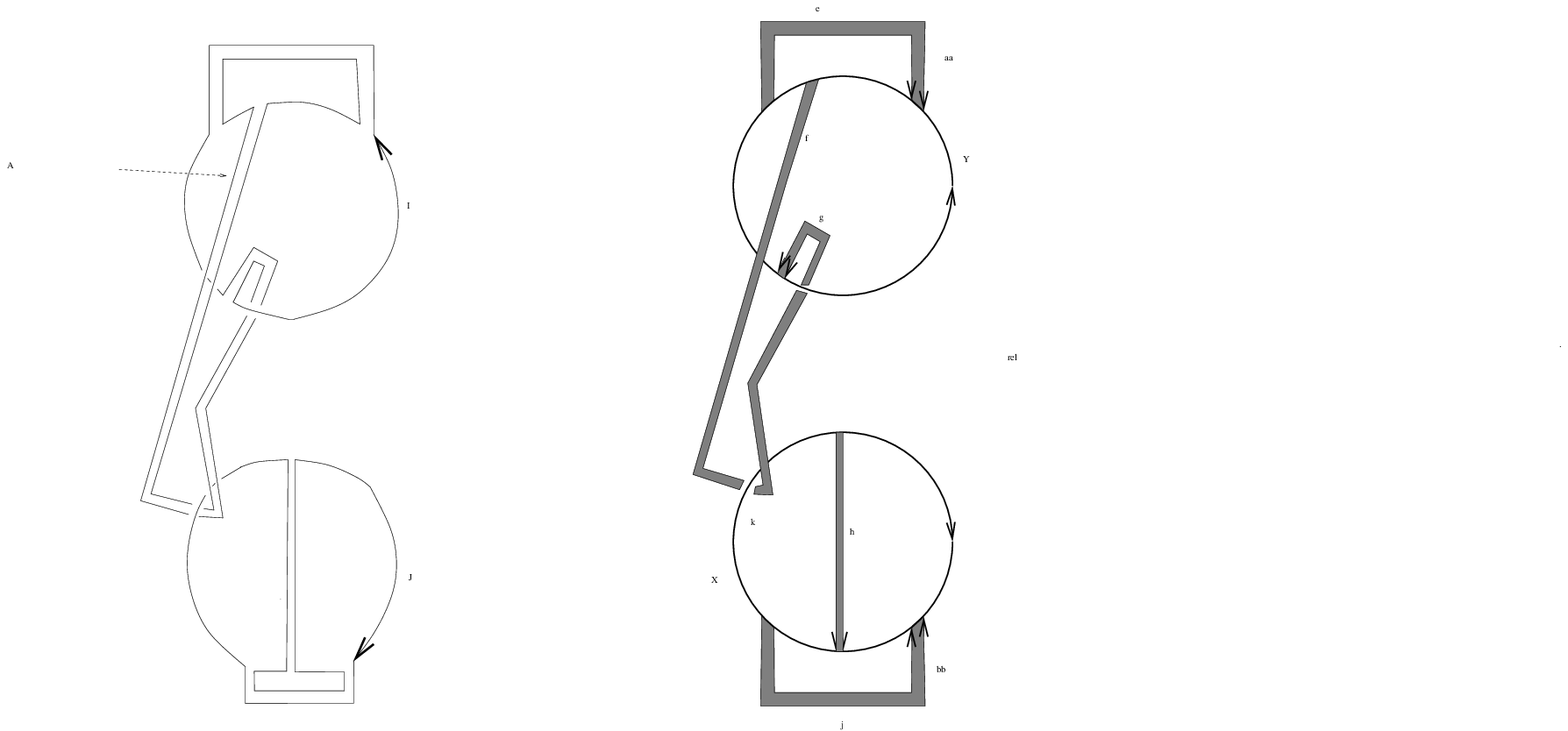}
\relabel{I}{$\s{Y}$}
\relabel{J}{$\s{X}$}
\relabel {e}{$\s{e}$}
\relabel {f}{$\s{f}$}
\relabel {g}{$\s{Y^{-1}X \t f}$}
\relabel {h}{$\s{h}$}
\relabel {j}{$\s{g}$}
\relabel {k}{$\s{X\t f}$} 
\relabel {X}{$\s{X}$}
\relabel {Y}{$\s{Y}$}
\relabel {A}{$\s{X^{-1}YX}$}
\relabel {rel}{{${\left \{ \begin{CD}  f(X \t f^{-1})  =1 \\ { (X \t f)  f^{-1}   }=1\\ \d(f)=1\\ \d(h)=1 \end{CD} \right . }$}}
\relabel {aa}{$ \d(e)=X^{-1}YXY^{-1}$}
\relabel {bb}{$ \d(g)=1$}
\endrelabelbox }
\caption{\label{Calc1a} Calculation of $\Pi_2(N,N^{(1)})$ where $N=S^4 \setminus \n(\S')$. Here $\S'$ is represented by the knot with bands $K'$ of figure \ref{movie2a}.}
\end{figure}

 Consider the knotted surface $\S'$ obtained from the  knot with bands $K'$ which appears in figure \ref{movie2a}.  Similarly with the Spun Hopf Link, this knotted surface is, topologically, diffeomorphic with the disjoint union of two {tori $S^1 \times S^1$}. 

Let us {determine} $\Pi_2(N,N^{(1)})$, where $N=S^4 \setminus \n(\S')$, {with the}  handle decomposition determined by $K'$. This calculation appears in figure \ref{Calc1a}. This permits us to conclude that  $\Pi_2(N,N^{(1)})$  is the {crossed module presented by the map $\{e,f,g,h\} \to F(X,Y)$, where {$f,h,g \mapsto 1$} and $e \mapsto X^{-1}YXY^{-1}$, considering the 2-relation $f(X \t f^{-1})=1.$}
 Recall that  $F(X,Y)$ is the free group on the variables $X$ and $Y$.

Let $\G=(G,E,\d,\t)$ be a finite crossed module. By using this calculation, {  we can determine} $I_\G\big(S^4\setminus \n(\S')\big)${. Suppose} that $\G=(G,E,\d,\t)$ is a finite crossed module with $G$ abelian and  $\d=1_G$. We have:
\begin{equation*}
I_\G(S^4\setminus {\n(\S')}){=\# G \# E \#\big \{X \in G; f \in E\left|  f(X \t f^{-1})  \right. =1_E \big \} } .
\end{equation*}
{This particularizes to {$I_\A(S^4 \setminus \n(\S'))=24$}, for the case when  $\A=(\Z_2,\Z_3,\d,\t)$ is the crossed module defined {in \ref{DDD}}.} This in particular proves that $\S'$ is knotted and that $\S'$ is not isotopic to  the Spun Hopf Link $\S$. 

 Note that the fundamental groups of {the complements of} $\S$ and $\S'$ each are isomorphic with $\Z\oplus \Z$ {(this can be inferred from  the presentations of the fundamental crossed modules of {them})}.  Therefore $\S$ and $\S'$ are two knotted surfaces (each a disjoint union of {two tori $S^1 \times S^1$)} with the same fundamental group of the complement, but distinguished by their crossed module invariants. This example appears in {\cite{FMK}.}

\section*{Acknowledgements}
This work had the financial support of  FCT (Portugal), post-doc grant number SFRH/BPD/17552/2004,  part of the research project POCTI/MAT/\\60352/2004 (''Quantum Topology''), also financed by FCT.

I would like to thank  a  referee of a previous version of this work; as well as Gustavo Granja,   Scott Carter and Roger Picken for useful comments.

\end{document}